%% file: paper.tex
\documentclass[12pt,reqno]{amsart}

\input next
\input epsf

\begin{document}

\bibliographystyle{plain}

\title[Patterns in permutations]{Enumeration of permutations containing
	a prescribed number of occurrences of a pattern of length 3}

\begin{abstract}
\input abstract
\end{abstract}

\author{Markus Fulmek}
\address{Institut f\"ur Mathematik der Universit\"at Wien\\
Strudlhofgasse 4, A-1090 Wien, Austria}
\email{{\tt Markus.Fulmek@Univie.Ac.At}\newline\leavevmode\indent
{\it WWW}: {\tt http://www.mat.univie.ac.at/\~{}mfulmek}
}
\date{\today}

\maketitle

\input macros

\section{Introduction}
\label{intro}
\input intro

\section{Results and conjectures}
\label{results}
\input results

\section{The proofs}
\label{proofs}
\input proof

\input proof1

\input proof2

\bibliography{paper}

\end{document}

%% file: abstract.tex
We consider the problem of enumerating the permutations
containing exactly $k$ occurrences of a pattern of length 3. This
enumeration has received a lot of interest recently, and there are a lot
of known results. This paper presents an alternative approach to
the problem, which yields a proof for a formula
which so far only was conjectured (by Noonan and Zeilberger). This approach
is based on bijections from permutations to certain lattice paths with
``jumps'', which were first considered by Krattenthaler.

%% file: macros.tex
\def\qp#1#2{[{#1}^{#2}]}
\def\sqp#1#2{s_{\qp{#1}{#2}}}
\def\pa#1{[#1]}
\def\sla#1{s_{\pa{#1}}}
\def\TP{P}
\def\TQ{\hat{P}}
\def\TS{\check{P}}
\def\minor#1#2#3#4#5{#1_{\{#2,#3\},\{#4,#5\}}}
\def\la{\lambda}
\def\si{\sigma}
\def\pprime{{\prime\prime}}

\def\q#1{[#1]_q}
\def\pq#1#2#3{(#1;#3)_{#2}}
\def\ceil#1{\lceil #1\rceil}
\def\floor#1{\lfloor #1\rfloor}
\newcommand{\HypFsimple}[2]{\sideset{_#1}{_#2}{F}}
\newcommand{\HypPhi}[6]{\sideset{_#1}{_#2}{\Phi}\!\left[\matrix  
#3\\#4\endmatrix;{\displaystyle #5,#6}\right]}

\def\bit{\begin{itemize}}
\def\eit{\end{itemize}}

\def\rlmum{left--to--right maximum}
\def\rlma{left--to--right maxima}
\def\jv{\mathcal J}
\def\tbase{$\tau$--base}

\newcommand{\SG}{{\mathcal S}}
\newcommand{\SGntr}[3]{\SG_{#1}\!\of{#2,\,#3}}
\newcommand{\occntr}[3]{s_{#1}\!\of{#2,\,#3}}
\newcommand{\GFxtr}[3]{F_{#2,\,#3}\!\of{#1}}

\let\epsilon\varepsilon
\newcommand{\esf}{{e}}                  
\newcommand{\hsf}{{h}}                  
\newcommand{\ksf}{{k}}                  
\newcommand{\schf}[2]{{s}_{#1}\of{#2,\bold{x}}}
\newcommand{\defeq}{:=}                                 
\newcommand{\tableau}{T}                                
\newcommand{\Ferrer}{B}                            
\newcommand{\digraph}{D}                                
\newcommand{\sdigraph}{S}                               
\newcommand{\horedge}{a}                                
\newcommand{\veredge}{\hat{a}}          
\newcommand{\fpath}{{\mathcal P}}                       
\newcommand{\nfpath}{{\cal N}}          
\newcommand{\nifpath}{{\cal P}_o}       
\newcommand{\GF}[1]{\bold{GF}\of{#1}}
\newcommand{\sgn}{\operatorname{sgn}}                   
\newcommand{\inv}{\operatorname{inv}}                   
\newcommand{\sschf}[2]{{sp}_{#1}\of{#2,\bold{x}}}
\newcommand{\isschf}[2]{{sp}_{n,m}\of{#1,#2;\bold{x},\bold z}}
\newcommand{\osschf}[2]{{sp}_{n,1}\of{#1,#2,\bold{x},z}}
\newcommand{\reflect}{{\bold{R}}}               
\newcommand{\mreflect}{\tilde{\bold{R}}}                
\newcommand{\oschf}[2]{{o}_{#1}\of{#2,\bold{x}}}
\newcommand{\oschftriv}[2]{{o}_{#1}\of{#2,\seqof{1,1,\dots,1}}}
\newcommand{\liff}{\text{ iff }}                
\newcommand{\limp}{\text{ implies }}
\renewcommand{\land}{\text{ and }}      

\newcommand{\thmref}[1]{Theorem~\ref{#1}}
\newcommand{\secref}[1]{Section~\ref{#1}}
\newcommand{\lemref}[1]{Lemma~\ref{#1}}
\newcommand{\figref}[1]{Figure~\ref{#1}}
\newcommand{\tabref}[1]{Table~\ref{#1}}

\newcommand{\seqof}[1]{\left\langle#1\right\rangle}
\newcommand{\setof}[1]{\left\{#1\right\}}
\newcommand{\of}[1]{\left(#1\right)}
\newcommand{\parof}[1]{\left(#1\right)}
\newcommand{\numof}[1]{\left|#1\right|}
\newcommand{\detof}[2]{\left|{#2}\right|_{#1\times #1}}
\newcommand{\firstcolumn}[1]{
\quad \vdots\quad }
\newcommand{\firstrow}[1]{
\quad \vdots\quad }
\newcommand{\brkof}[1]{\left[#1\right]}
\newcommand{\intof}[1]{\lfloor #1\rfloor}
\newcommand{\absof}[1]{\left|#1\right|}
\newcommand{\sgnof}[1]{\operatorname{sgn}\of{#1}}

\newtheorem{thm}{Theorem}
\newtheorem{lem}[thm]{Lemma}
\newtheorem{cor}[thm]{Corollary}
\newtheorem{dfn}[thm]{Definition}
\newtheorem{obs}[thm]{Observation}
\newtheorem{rem}[thm]{Remark}
\newtheorem{con}[thm]{Conjecture}

\setlength{\unitlength}{0.6cm}
\newsavebox{\schraffenbox}
\savebox{\schraffenbox}(1,1)[lb]{
\put(0,0){\line(1,1){1}}
\put(0,0.5){\line(1,1){0.5}}
\put(0.5,0){\line(1,1){0.5}}
}
\newsavebox{\Dyckbox}
\savebox{\Dyckbox}(2,1)[lb]{
\put(0,0){\line(1,1){1}}
\put(1,1){\line(1,-1){1}}
\put(2,0){\line(-1,0){2}}
\put(0.8,0.2){{\bf c}}
}
\newsavebox{\upstepbox}
\savebox{\upstepbox}(1,1)[lb]{
\linethickness{2pt}
\put(0,0){\line(1,1){1}}
\put(0,0){\circle*{0.2}}
\put(1,1){\circle*{0.2}}
\put(-0.05,0.65){$\scriptstyle\sqrt x$}
}
\def\dpath#1{$\dosteplist#1\mendlist$}
\def\mendlist{\mendlist}
\def\dosteplist{\afterassignment\handlenextstep\let\next=}
\def\handlenextstep{
\ifx\next\mendlist
	\let\next=\relax
\else
	\ifx\next u
   	/
  	\else
   	\ifx\next j
    		\vert
   	\else
    		\ifx\next d
     			\backslash
			\else
				\errmessage{Wrong symbol?}
     		\fi
		\fi
	\fi
	\let\next=\dosteplist
\fi
\next
}

\def\useg#1{$\underbrace{/\dots/}_{#1}$}
\def\dseg#1{$\underbrace{\backslash\dots\backslash}_{#1}$}
\def\jseg#1{$\underbrace{\vert\dots\vert}_{#1}$}
\def\usega{$/\dots/$}
\def\dsega{$\backslash\dots\backslash$}
\def\jsega{$\vert\dots\vert$}
\def\dstepat#1{$\backslash^{\!\!#1}$}
\def\dsegat#1#2{$\backslash\cdot(#1)\cdot$\dstepat{#2}}
\def\dsegata#1{$\backslash\dots$\dstepat{#1}}

%% file: intro.tex
\subsection{Patterns in permutations}
We consider the group $\SG_n$ of permutations of
the set $\{1,\dots,n\}$.  Let $\tau$ be an arbitrary permutation in
$\SG_k$ ($k \geq 2 $, in order to avoid trivial cases).
We say that some permutation
$\rho=\of{\rho_1,\dots,\rho_n}\in\SG_n$ {\em contains\/} the {\em pattern\/} $\tau$,
if there exists a subword $\sigma=\of{\rho_{i_1},\dots,\rho_{i_k}}$ in $\rho$, such that the entries of the subword $\sigma$ appear in the same
``relative order'' as the entries of the permutation $\tau$.
More formally stated, the subword $\sigma$ can be transformed
into the permutation $\tau$ by the following construction. Replace the smallest
element in $\sigma$ by 1, the second--smallest by 2, the third--smallest
by 3, and so on.

If  $\rho\in\SG_n$ does not contain pattern $\tau$,
we say that $\rho$ is {\em $\tau$--avoiding\/}.

For example, the permutation $\rho=(1,5,2,4,3)$ contains the
pattern $\tau=(3,1,2)$,
since the subword $\sigma=(\rho_2,\rho_3,\rho_4)=(5,2,4)$ shows the same
relative ordering of elements as $\tau$. First comes the largest element
(i.e., 5), then the  smallest (i.e., 2), and then the second--largest
(i.e., 4).

\begin{figure}
\caption{Graphical illustration. Permutation $(1,5,2,4,3)$ contains
	exactly 2 $(3,1,2)$--patterns.}
\label{fig:patterns}
\begin{picture}(6,6)(0,0)
\put(0.5,5.5){\line(1,0){5}}
\put(0.5,4.5){\line(1,0){5}}
\put(0.5,3.5){\line(1,0){5}}
\put(0.5,2.5){\line(1,0){5}}
\put(0.5,1.5){\line(1,0){5}}
\put(0.5,0.5){\line(1,0){5}}
\put(1,0){{\small 1}}
\put(2,0){{\small 2}}
\put(3,0){{\small 3}}
\put(4,0){{\small 4}}
\put(5,0){{\small 5}}
\put(5.5,0.5){\line(0,1){5}}
\put(4.5,0.5){\line(0,1){5}}
\put(3.5,0.5){\line(0,1){5}}
\put(2.5,0.5){\line(0,1){5}}
\put(1.5,0.5){\line(0,1){5}}
\put(0.5,0.5){\line(0,1){5}}
\put(0,1){{\small 1}}
\put(0,5){{\small 5}}
\put(0,3){{\small 3}}
\put(0,4){{\small 4}}
\put(0,2){{\small 2}}
\put(1,1){\circle*{0.3}}
\put(2,5){\circle*{0.3}}
\put(3,2){\circle*{0.3}}
\put(4,4){\circle*{0.3}}
\put(5,3){\circle*{0.3}}
\put(1.9,4.7){\line(1,-3){1.1}}
\put(2.4,4.8){\line(1,-3){0.7}}
\put(3.1,2.7){\vector(1,2){0.7}}
\put(3.0,1.4){\vector(1,1){1.7}}
\end{picture}

\end{figure}
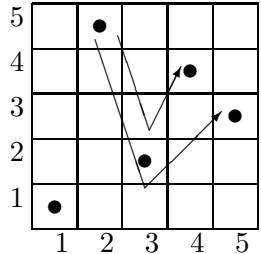

The enumeration of the set $\SGntr{n}{\tau}{r}$ of permutations from
$\SG_n$, which contain the fixed pattern $\tau$ precisely $r$ times,
has recently received considerable
interest \cite{barcucci:perm,barcucci:cn,barcucci:incr,bona:1342,bona:precursive,bona:length4,bona:132,egge:fibonacci,krat:dyck,mansour-vainshtein:continued,mansour-vainshtein:counting,mansour-vainshtein:restricted,mansour-vainshtein:restricted1,mansour-vainshtein:layered,noonan:one123,noonan-zeilberger:enum123,robertson:length3,robertson-wilf-zeilberger:patterns,simion-schmidt:restricted,west:phd}.

In our above example $\rho=(1,5,2,4,3)$, the pattern $(3,1,2)$ occurs
exactly twice,
since the subword $\sigma=(\rho_2,\rho_3,\rho_5)=(5,2,3)$, too, shows the same
relative ordering as $(3,1,2)$. Hence $(1,5,2,4,3)$ belongs to
$\SGntr{5}{(3,1,2)}{2}$.

We shall denote the cardinality of
$\SGntr{n}{\tau}{r}$ by
$\occntr{n}{\tau}{r}$, where we set, by convention,
$\occntr{0}{\tau}{0}=\occntr{1}{\tau}{0}=1$ and
$\occntr{0}{\tau}{r}=\occntr{1}{\tau}{r}=0$ for $r>0$ (recall that $\tau$
is of length at least 2). We denote the generating function of these numbers
by
$
\GFxtr{x}{\tau}{r} := \sum_{n=0}^\infty \occntr{n}{\tau}{r} x^n
$

\subsection{Outline of this paper}

In this paper, we consider patterns of length 3. It is well known
that there are only two patterns in $\SG_3$ which are essentially different,
namely $(3,2,1)$ and $(3,1,2)$. (We shall give a simple argument for this
fact right below). For both of these patterns, we make use of a bijective
construction which basically ``translates'' permutations into certain
``generalized'' Dyck--paths, where ``jumps'' are allowed.
This construction was first considered by
Krattenthaler \cite[Lemma $\phi$ and Lemma $\Psi$]{krat:dyck}.
Using the well--known generating function for Dyck--paths, we are able to
obtain the generating functions $\GFxtr{x}{(3,1,2)}{r}$ and $\GFxtr{x}{(3,2,1)}{r}$ for $r=0,1,2$ in a uniform way. 
Especially, for $r=0$ we thus obtain a very simple ``graphical'' proof of
the formulas for $\occntr{n}{(3,1,2)}{0}$ and
$\occntr{n}{(3,2,1)}{0}$ (which, of course, are very well known, see
\cite[Exercises 6.19.ee and 6.19.ff]{stanley:ec2}).
The derivation
of $\GFxtr{x}{(3,2,1)}{2}$ seems to be new (it was conjectured by
Noonan and Zeilberger \cite{noonan-zeilberger:enum123}).

In \secref{results} we present the respective formulas. In \secref{proofs}
we introduce the bijective construction and show how the generating
functions can be derived. 

\subsection{Graphical approach to permutation patterns}
``Graphical arguments'' play a key role in our presentation. We introduce
the {\em permutation graph\/} $G\of\rho$ (see \cite[p.~71]{stanley:ec1})
as a simple way for visualizing some permutation $\rho\in\SG_n$:
\begin{equation*}
G\of\rho:=
	\{(1,\rho_1),\dots,(n,\rho_n)\}\subseteq\{1,\dots,n\}\times\{1,\dots,n\}.
\end{equation*}
Figure~\ref{fig:patterns} illustrates this concept. The set 
$\{1,\dots,n\}\times\{1,\dots,n\}$ is visualized by a rectangular grid of
square cells, the elements $(i,\rho_i)$
of the graph are indicated by dots in the respective cell. Note that the two $(3,1,2)$--patterns in $(1,5,2,4,3)$
appear as ``hook--like configurations'', as indicated by the arrows
in the picture.

Clearly, $G\of\rho$ is just another way to denote the permutation $\rho$.
Hence, in the following we shall not make much difference between a
permutation and its graph.

We call the ``horizontal
line'' $\setof{(1,i),\dots,(n,i)}$ the $i$--th row of the graph, and
the ``vertical line'' $\setof{(j,1),\dots,(j,n)}$ the $j$--th column of
the graph. Note that the rows of $G\of\rho$ correspond to the
rows of the {\em permutation matrix\/} of $\rho$
in reverse order. Clearly, each
row and each column of the graph of some permutations contains precisely
one dot. 

\begin{figure}
\caption{Graphical illustration of the orbit of $(3,1,2)$ under rotations
of the permutation graph}
\label{fig:312orbit}
\begin{picture}(10,10)(0,0)
\put(0.5,5.5) {
\begin{picture}(4,4)(0,0)
\put(0.5,3.5){\line(1,0){3}}
\put(0.5,2.5){\line(1,0){3}}
\put(0.5,1.5){\line(1,0){3}}
\put(0.5,0.5){\line(1,0){3}}
\put(1,0){{\small 1}}
\put(2,0){{\small 2}}
\put(3,0){{\small 3}}
\put(3.5,0.5){\line(0,1){3}}
\put(2.5,0.5){\line(0,1){3}}
\put(1.5,0.5){\line(0,1){3}}
\put(0.5,0.5){\line(0,1){3}}
\put(0,3){{\small 3}}
\put(0,1){{\small 1}}
\put(0,2){{\small 2}}
\put(1,3){\circle*{0.3}}
\put(2,1){\circle*{0.3}}
\put(3,2){\circle*{0.3}}
\end{picture}
}
\put(5.5,5.5) {
\begin{picture}(4,4)(0,0)
\put(0.5,3.5){\line(1,0){3}}
\put(0.5,2.5){\line(1,0){3}}
\put(0.5,1.5){\line(1,0){3}}
\put(0.5,0.5){\line(1,0){3}}
\put(1,0){{\small 1}}
\put(2,0){{\small 2}}
\put(3,0){{\small 3}}
\put(3.5,0.5){\line(0,1){3}}
\put(2.5,0.5){\line(0,1){3}}
\put(1.5,0.5){\line(0,1){3}}
\put(0.5,0.5){\line(0,1){3}}
\put(0,1){{\small 1}}
\put(0,3){{\small 3}}
\put(0,2){{\small 2}}
\put(1,1){\circle*{0.3}}
\put(2,3){\circle*{0.3}}
\put(3,2){\circle*{0.3}}
\end{picture}
}
\put(0.5,0.5) {
\begin{picture}(4,4)(0,0)
\put(0.5,3.5){\line(1,0){3}}
\put(0.5,2.5){\line(1,0){3}}
\put(0.5,1.5){\line(1,0){3}}
\put(0.5,0.5){\line(1,0){3}}
\put(1,0){{\small 1}}
\put(2,0){{\small 2}}
\put(3,0){{\small 3}}
\put(3.5,0.5){\line(0,1){3}}
\put(2.5,0.5){\line(0,1){3}}
\put(1.5,0.5){\line(0,1){3}}
\put(0.5,0.5){\line(0,1){3}}
\put(0,2){{\small 2}}
\put(0,3){{\small 3}}
\put(0,1){{\small 1}}
\put(1,2){\circle*{0.3}}
\put(2,3){\circle*{0.3}}
\put(3,1){\circle*{0.3}}
\end{picture}
}
\put(5.5,0.5) {
\begin{picture}(4,4)(0,0)
\put(0.5,3.5){\line(1,0){3}}
\put(0.5,2.5){\line(1,0){3}}
\put(0.5,1.5){\line(1,0){3}}
\put(0.5,0.5){\line(1,0){3}}
\put(1,0){{\small 1}}
\put(2,0){{\small 2}}
\put(3,0){{\small 3}}
\put(3.5,0.5){\line(0,1){3}}
\put(2.5,0.5){\line(0,1){3}}
\put(1.5,0.5){\line(0,1){3}}
\put(0.5,0.5){\line(0,1){3}}
\put(0,2){{\small 2}}
\put(0,1){{\small 1}}
\put(0,3){{\small 3}}
\put(1,2){\circle*{0.3}}
\put(2,1){\circle*{0.3}}
\put(3,3){\circle*{0.3}}
\end{picture}
}
\end{picture}
\end{figure}
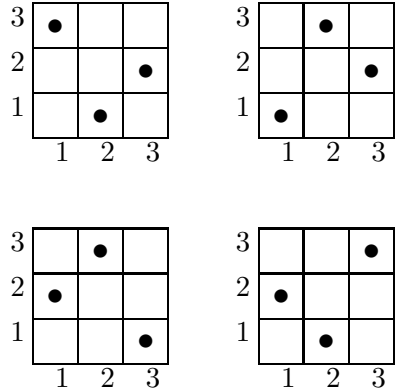

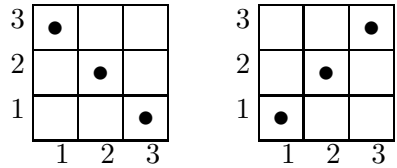
\begin{figure}
\caption{Graphical illustration of the orbit of $(3,2,1)$ under rotations
of the permutation graph}
\label{fig:321orbit}
\begin{picture}(10,5)(0,0)
\put(0.5,0.5) {
\begin{picture}(4,4)(0,0)
\put(0.5,3.5){\line(1,0){3}}
\put(0.5,2.5){\line(1,0){3}}
\put(0.5,1.5){\line(1,0){3}}
\put(0.5,0.5){\line(1,0){3}}
\put(1,0){{\small 1}}
\put(2,0){{\small 2}}
\put(3,0){{\small 3}}
\put(3.5,0.5){\line(0,1){3}}
\put(2.5,0.5){\line(0,1){3}}
\put(1.5,0.5){\line(0,1){3}}
\put(0.5,0.5){\line(0,1){3}}
\put(0,3){{\small 3}}
\put(0,2){{\small 2}}
\put(0,1){{\small 1}}
\put(1,3){\circle*{0.3}}
\put(2,2){\circle*{0.3}}
\put(3,1){\circle*{0.3}}
\end{picture}
}
\put(5.5,0.5) {
\begin{picture}(4,4)(0,0)
\put(0.5,3.5){\line(1,0){3}}
\put(0.5,2.5){\line(1,0){3}}
\put(0.5,1.5){\line(1,0){3}}
\put(0.5,0.5){\line(1,0){3}}
\put(1,0){{\small 1}}
\put(2,0){{\small 2}}
\put(3,0){{\small 3}}
\put(3.5,0.5){\line(0,1){3}}
\put(2.5,0.5){\line(0,1){3}}
\put(1.5,0.5){\line(0,1){3}}
\put(0.5,0.5){\line(0,1){3}}
\put(0,1){{\small 1}}
\put(0,2){{\small 2}}
\put(0,3){{\small 3}}
\put(1,1){\circle*{0.3}}
\put(2,2){\circle*{0.3}}
\put(3,3){\circle*{0.3}}
\end{picture}
}
\end{picture}
\end{figure}

A simple geometric argument now shows that, when talking about
patterns of length 3, we may restrict our attention to 2 ``essentially
different'' patterns,
namely $(3,2,1)$ and $(3,1,2)$. Given the permutation graph
of an arbitrary permutation $\tau$, consider the rotations by $0,\frac{\pi}{2},
\pi$ and
$\frac{3\pi}{2}$. Clearly, every pattern contained in $\tau$ is
rotated accordingly, so the numbers $\occntr{n}{\tau}{r}$ are the same as
the numbers $\occntr{n}{\tau^\prime}{r}$, where $\tau^\prime$ denotes the permutation
obtained by rotating the graph of $\tau$. Therefore,
we only have to consider the different {\em orbits\/} of the permutation graphs
of $\SG_3$ under the action of this rotation group.
It is easy to see that there are precisely two of them,
one containing $(3,1,2)$ (see Figure~\ref{fig:312orbit}),
the other containing $(3,2,1)$ (see Figure~\ref{fig:321orbit}).

%% file: results.tex
\subsection{Formulas for 0,1, and 2 occurrences of pattern $\tau$}
In the following, we summarize the formulas and generating functions we
are going to prove later. We start with a very well known fact:
\begin{equation}
\occntr{n}{(3,1,2)}{0} = \occntr{n}{(3,2,1)}{0} = \frac{1}{n+1}\binom{2n}{n}. \label{eq:s-0}
\end{equation}
These are the Catalan numbers, which go like this: $1$, $1$, $2$, $5$, $14$, $42$, $132$, $429$, $1430$, $4862$, $16796$, $58786$,
$208012$, $\dots$.
The corresponding generating function is
\begin{equation}
\GFxtr{x}{(3,1,2)}{0} = \GFxtr{x}{(3,2,1)}{0} = \frac{1}{2x}-\frac{1}{2x}\sqrt{1-4x}.
	\label{eq:S-0}
\end{equation}

The next formula was first proved by B\'ona \cite{bona:132}: 
\begin{equation}
\occntr{n}{(3,1,2)}{1} = \binom{2n-3}{n-3}. \label{eq:s312-1}
\end{equation}
The numbers go like this: $0$, $0$, $0$, $1$, $5$, $21$, $84$, $330$, $1287$, $5005$, $19448$, $75582$,
$293930$, $\dots$.
The corresponding generating function is
\begin{equation}
\GFxtr{x}{(3,1,2)}{1} = \frac{x-1}{2}-\frac{3x-1}{2}(1-4x)^{-1/2}.
	\label{eq:S312-1}
\end{equation}

The next formula 
was first
proved by Mansour and Vainshtein \cite{mansour-vainshtein:counting}:
\begin{equation}
\occntr{n}{(3,1,2)}{2} = \binom{2n-6}{n-4}\frac{n^3+17n^2-80n+80}{2n(n-1)}, 
	\label{eq:s312-2}
\end{equation}
The
numbers go like $0$, $0$, $0$, $0$, $4$, $23$, $107$, $464$, $1950$, $8063$, $33033$, $134576$, $546312$, $\dots$.
The corresponding generating function is
\begin{multline}
\GFxtr{x}{(3,1,2)}{2} = \frac{x^2+3x-2}{2}\\
	+\frac{2x^4-4x^3+29x^2-15x+2}{2}(1-4x)^{-3/2}. 
	\label{eq:S312-2}
\end{multline}

\begin{rem}
Mansour and Vainshtein
\cite[Corollary 3]{mansour-vainshtein:counting}
are able to compute the numbers $\occntr{n}{(3,1,2)}{r}$
and generating functions $\GFxtr{x}{(3,1,2)}{r}$ for $r$ up to 6.
\end{rem}

The next formula was first
proved by Noonan \cite{noonan:one123}:
\begin{equation}
\occntr{n}{(3,2,1)}{1} = \frac{3}{n}\binom{2n}{n-3}. \label{eq:s321-1}
\end{equation}
These numbers go like $0$, $0$, $0$, $1$, $6$, $27$, $110$, $429$, $1638$, $6188$, $23256$, $87210,
326876$, $\dots$.
The corresponding generating function is
\begin{equation}
\GFxtr{x}{(3,2,1)}{1} = -\frac{2x^3-9x^2+6x-1}{2x^3}-\frac{\sqrt{1-4x}(3x^2-4x+1)}{2x^3}.
	\label{eq:S321-1}
\end{equation}

The next formula was conjectured by Noonan and Zeilberger
\cite{noonan-zeilberger:enum123}; to the best of my knowledge the first
proof is contained in this paper:
\begin{equation}
\occntr{n}{(3,2,1)}{2} = \frac{59n^2+117n+100}{2n(2n-1)(n+5)}\binom{2n}{n-4}.
	\label{eq:s321-2}
\end{equation}
These numbers go like $0$, $0$, $0$, $0$, $3$, $24$, $133$, $635$, $2807$, $11864$, $48756$, $196707,
783750$, $\dots$.
The corresponding generating function is
\begin{multline}
\GFxtr{x}{(3,2,1)}{2} = -\frac{5x^5-7x^4+17x^3-20x^2+8x-1}{2x^5}\\+
	\frac{\sqrt{1-4x}(x^5-3x^4+5x^3-10x^2+6x-1)}{2x^5}.
	\label{eq:S321-2}
\end{multline}

\begin{rem}
Toufik Mansour \cite{mansour:private} informed me that he found
a proof of formula~\ref{eq:s321-2}, too.
\end{rem}

\subsection{Further conjectures}

We state the following conjectures concerning generating functions.
\begin{con}
The generating function $\GFxtr{x}{(3,2,1)}{3}$ is 
\begin{multline}
\GFxtr{x}{(3,2,1)}{3} = -\frac{(-1 + 10x - 33x^2 + 32x^3 + 31x^4 - 70x^5 + 
   35x^6 - 2x^8)}{2x^7}\\+
	\frac{\sqrt{1-4x}(-1 + 8x - 19x^2 + 6x^3 + 27x^4 - 28x^5 + 7x^6 + 
    2x^7)}{2x^7}.
	\label{eq:S321-3}
\end{multline}
\end{con}

\begin{con}
The generating function $\GFxtr{x}{(3,2,1)}{4}$ is 
\begin{multline}
\GFxtr{x}{(3,2,1)}{4} = \frac{-1}{2x^9}\bigl[
	(-1 + 12x - 50x^2 + 65x^3 + 107x^4 - 437x^5 + 588x^6\\
		 - 492x^7 + 314x^8 - 108x^9 + 3x^{10})\\+
	\sqrt{1-4x}(-1 + 10x - 32x^2 + 17x^3 + 107x^4 - 245x^5 + 256x^6 - 
    192x^7 + 102x^8\\
	 	- 18x^9 - x^{10})
	\bigr].
	\label{eq:S321-4}
\end{multline}
\end{con}

\subsection{General form of the formulas}

By inspection of the formulas listed above, one is immediately led
to conjecture a simple general form for the generating functions
$\GFxtr{x}{(3,1,2)}{r}$ and $\GFxtr{x}{(3,2,1)}{2}$.
The ``$(3,1,2)$--part'' was proved by
B\'ona~\cite[Proposition~1 and Lemma~3]{bona:precursive}.

\begin{thm}[B\'ona]
The generating functions $\GFxtr{x}{(3,1,2)}{r}$ are rational functions
in $x$ and
$\sqrt{1-4x}$. Moreover, when written in smallest terms, the denominator
of $\GFxtr{x}{(3,1,2)}{r}$ is equal to $\of{\sqrt{1-4x}}^{2r-1}$.
\end{thm}

For the ``$(3,2,1)$--part'', we state the following conjecture:

\begin{con}
The generating functions $\GFxtr{x}{(3,2,1)}{r}$ is of the form
\begin{equation*}
\frac{1}{2x^{2r+1}}\of{P_r\of x + \sqrt{1-4x}\, Q_r\of x},
\end{equation*}
where $P_r$ and $Q_r$ are polynomials.
\end{con}

\begin{rem}
Toufik Mansour \cite{mansour:private} informed me that he found
a proof for the general form of the generating function
$\GFxtr{x}{(3,2,1)}{r}$.
\end{rem}

%% file: proof.tex
\subsection{The case of pattern--avoidance}
\label{sec:occ0}
Let us start with a simple proof of \eqref{eq:s-0}. Recall that the
Catalan numbers $C_n$ arise in the enumeration of Dyck paths, which are
paths in the integer lattice consisting of ``up--steps'' $(1,1)$ and ``down--steps'' $(1,-1)$, which start in $(0,0)$ and end at $(2n,0)$,
such that the path never goes below the horizontal axis. More precisely, we have
\begin{equation*}
C_{n} = \frac{1}{n+1}\binom{2n}{n} = \text{\# of Dyck--paths from $(0,0)$ to $(2n,0)$}.
\end{equation*}

So \eqref{eq:s-0} will follow, if we can establish a bijection between
permutations in $\SG_n$ which do avoid the pattern $(3,2,1)$
or $(3,1,2)$, respectively, and Dyck--paths from $(0,0)$ to $(2n,0)$. Denote
the set of such paths by $D_n$.

\bigskip
In the following, $\tau$ will always denote one of the patterns $\of{3,1,2}$
or $\of{3,2,1}$. If some assertion A is valid for both patterns, we shall
simply say ``A is valid for $\tau$'' instead of ``A is valid for $\of{3,1,2}$
and for $\of{3,2,1}$''. For an entry $\rho_i$ in some permutation 
$\rho=\of{\rho_1,\dots,\rho_n}\in\SG_n$, we call the entries $\rho_j$
with $j<i$ the entries {\em to the left\/}, and with $j\leq i$ the
entries {\em weakly to the left\/}, and likewise for entries (weakly)
to the right.

\begin{dfn}
Let $\rho=\of{\rho_1,\dots,\rho_n}\in\SG_n$ be an arbitrary permutation.
An entry $\rho_i$ is called a {\em \rlmum\/}, if $\rho_i$ is greater than
all entries to the left of it; i.e., if $\rho_i>\rho_j$ for all $j < i$.
Entries which are {\em not\/} \rlma\ are called {\em remaining entries\/}.
\end{dfn}


\begin{obs}
By definition, for a subword $\of{\sigma_1,\sigma_2,\sigma_3}$ in $\rho$
which establishes an occurrence of the pattern $\tau$, neither entry $\sigma_2$ nor $\sigma_3$
can be a \rlmum\ in $\rho$. On the other hand, the entry $\sigma_1$ is either a
\rlmum\ itself, or there is a \rlmum\ $\rho_k$ to the left of $\sigma_1$,
such that $\of{\rho_k,\sigma_2,\sigma_3}$ gives another $\tau$--occurrence.
In fact, this is the case for {\em all\/} \rlma\ $\rho_k$ weakly to the left
of $\sigma_2$, such that $\rho_k>\max\of{\sigma_2,\sigma_3}$.
\end{obs}
To establish our desired bijections, we first consider a single mapping
\begin{equation*}
\psi:\;\SG_n\rightarrow D_n,
\end{equation*}
which will turn out to yield both bijections for $\tau=(3,2,1)$ and
$\tau=(3,1,2)$ by simply restricting to  $\SGntr{n}{\tau}{0}$.
$\psi$ is defined
by the following ``graphical'' construction: For an arbitrary permutation
$\rho\in\SG_n$,
\bit  
\item draw the permutation graph of $\rho$ and mark the \rlma\ (i.e.,
	the points $(i,\rho_i)$ with $\rho_i>\rho_j$ for all $j < i$),
\item consider the points $(x,y)\in\{1,\dots,n\}\times\{1,\dots,n\}$,
	which lie below and to the right of some \rlmum\ (i.e., $x\geq i$ and
	$y\leq \rho_i$ for some \rlmum\ $(i,\rho_i)$),
\item rotate the ``polygonal figure'' formed by these
	points by $-\frac{\pi}{4}$ and consider its upper boundary, which appears
	as some path $p$ from $(0,0)$ to $(2n,0)$,
\item set $\psi\of{\rho} = p$.
\eit

Figure~\ref{fig:bijection0} gives an illustration of this construction.
The \rlma\ are marked by circles, the boundary line of
the polygonal figure is shown as bold line, and the corresponding path is
drawn below. Note that the entries of the permutation correspond to
the down--steps of the path (this correspondence is indicated by the
labels in Figure~\ref{fig:bijection0}); the \rlma\ of the permutation
correspond to the ``peaks'' of the path. In the following, we will more
or less identify the $i$--th entry of $\rho$ and the $i$--th down--step of
$\psi\of\rho$.

\input fig1

We have to show that the resulting path $p=\psi\of\rho$ is
actually in $D_n$, i.e., $p$ does not go below the horizontal axis. For
convenience, we introduce some notation concerning paths.

\begin{dfn}
\label{dfn:prec-rlmum}
A down--step from $\of{x,h+1}$ to $\of{x+1,h}$ is said to have (or be of)
height $h$. By abuse of notation, we will also say that the entry of the
permutation corresponding to such down--step has height $h$.

The $i$--th down--step in $\psi\of\rho$
corresponds to the $i$--th entry of $\rho$. There must exist \rlma\ $\rho_j$
weakly to the left of this entry (i.e., $j\leq i$). We call the
right--most of these the {\em preceding \rlmum\/} with respect
to the down--step; it corresponds to the peak  weakly to the left
of the down--step.
\end{dfn}

So what we have to show is that all down--steps have nonnegative heights.
First note that the height of a \rlmum\ $\rho_m$ is exactly the number of entries to the right which are smaller; i.e.,
\begin{equation*}
\#\{j: m < j \leq n\text{ and } \rho_j<\rho_m\},
\end{equation*}
which clearly is nonnegative.
More generally, let $h_k$ be the height of the $k$--th down--step, and let
$\rho_m$ be the preceding \rlmum\ of height $h_m$. It is easy to see that
$h_k = h_m - (k-m)$, which is nonnegative, too, since $(k-m)$ cannot exceed $h_m$.
Hence $p=\psi\of\rho$ is indeed in $D_n$ for all $\rho\in\SG_n$.

Next we show that the mapping $\psi$, restricted to $\SGntr{n}{\tau}{0}$,
is a bijection, for $\tau=(3,2,1)$ as well as for $\tau=(3,1,2)$. Thus \eqref{eq:s-0} will follow immediately. Given an
arbitrary Dyck--path $p$ in $D_n$, we will construct a permutation
$\rho\in\SG_n$ such that
\bit
\item $\psi\of\rho=p$,
\item $\rho$ is the {\em unique\/} $\tau$--avoiding permutation in
	$\psi^{-1}\of p$.
\eit

 We start by rotating path $p$ by $\frac{\pi}{4}$, thus obtaining a
``configuration of \rlma'', which we may view as an ``incomplete
permutation graph''.

While there are empty rows and columns in this ``incomplete graph'',
repeat the following step for the case $\tau=(3,2,1)$:
\begin{quote}
Consider the {\em left--most\/} column $j$, which does not yet contain a dot, and the {\em lowest\/} row $i$,
which does not yet contain a dot, and put a dot into the cell $(i,j)$.
\end{quote}
Observe that a permutation obtained
from this construction is $(3,2,1)$-avoiding.
Figure~\ref{fig:bijection0-1} illustrates this construction for the
``configuration of \rlma'' from Figure~\ref{fig:bijection0}.

For the case $\tau=(3,1,2)$, repeat the following step:
\begin{quote}
Consider the left--most column $j$, which does not yet contain a dot, and the {\em highest\/} row $i$, which does not yet contain a dot {\em and\/}
which lies {\em below\/} the preceding \rlmum, and put a dot into the cell $(i,j)$.
\end{quote}
Observe that a permutation obtained
from this construction is $(3,1,2)$-avoiding,
as is illustrated by Figure~\ref{fig:bijection0}

These constructions will give completed graphs (and thus some
permutation $\rho$) in any case; however, in order to have $\psi\of\rho=p$,
there must hold the following condition:
\begin{quote}
Whenever there is a column $j$ which does not yet contain a dot, then there is
also some row $i$ without a dot, which lies {\em below\/} the
\rlmum\ {\em preceding\/} (in the same sense as before) column $j$.
\end{quote}
But this is precisely the condition, that the path
$p$ does not go below the horizontal axis.

\input fig2

The fact that the permutation $\rho$ is indeed unique is best seen ``by
inspection'', see Figures~\ref{fig:bijection0} (for $(3,1,2)$--avoiding
permutations) and \ref{fig:bijection0-1} (for $(3,2,1)$--avoiding
permutations).
Consider an arbitrary ``configuration of \rlma'' and
observe:
\begin{itemize}
\item If a permutation is $(3,1,2)$--avoiding, the ``remaining elements''
	must be inserted in a ``descending'' way,

\item if a permutation is $(3,2,1)$--avoiding, the ``remaining elements''
	must be inserted in an ``ascending'' way.\qed
\end{itemize}

\subsection{Krattenthaler's bijections}
\label{sec:kratt}

The bijective proof for the case of zero occurrences of length--$3$--patterns
was very simple. While the formal description of the constructions took
some pages of text, the main idea is easily obtained by mere ``inspection'' of
Figure~\ref{fig:bijection0} and Figure~\ref{fig:bijection0-1}.
This ``graphical view'' will turn out to be fruitful also for the
cases of one and more occurrences of patterns; however, we do not have
such a simple and uniform ``graphical'' definition any more.
We now have to distinguish between the
case of $(3,2,1)$--patterns and $(3,1,2)$--patterns. The respective
constructions were given by Krattenthaler \cite{krat:dyck}.

\begin{dfn}
\label{dfn:jumppaths}
Consider the set of lattice paths from $(0,0)$ to $(2n+s,0)$, consisting of $(n+s)$
up--steps $(1,1)$, $n$ down--steps $(1,-1)$, and $s$ down--jumps $(0,-1)$,
which do never go below the horizontal axis. Denote this set
of ``generalized Dyck--paths with down--jumps'' by $D_{n,s}$,
denote the union $\bigcup_{s=0}^{\infty}D_{n,s}$ by $D^\prime_n$
and the union $\bigcup_{n=0}^{\infty}D^\prime_n$ by $D^\prime$.

A maximal seqence of $d$ consecutive down--jumps is called a
{\em jump of depth $d$\/}.

\end{dfn}

Such paths with jumps will be the main object of the following considerations.
We will refer to them simply as ``lattice paths'' or (even simpler) as
``paths''. In order to avoid complicated verbal descriptions of such paths
(or segments of
paths), we introduce a ``graphical notation'' by denoting an up--step
by 
\epsffile{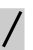}, a down--step by 
\epsffile{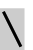}, and a down--jump by 
\epsffile{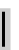}. For illustration, see Figure~\ref{fig:bijections}:
The two Dyck--paths with jumps shown there read
\epsffile{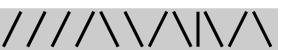}\ and
\epsffile{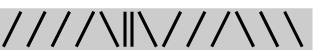}, respectively, in this ``graphical notation''.

For specifying the length of some sequence of consecutive steps, we
introduce a notation by means of example:
\epsffile{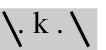}\ denotes a sequence of $k$ downsteps.
A down--step of height $h$ (in the sense of Definition~\ref{dfn:prec-rlmum})
is denoted by
\epsffile{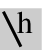}.

\begin{dfn}
\label{dfn:psi312}
Define the mapping
$$\psi_{(3,1,2)}:\;\SG_n\rightarrow D^\prime_n$$
by the following construction: Read $\rho=\of{\rho_1,\dots,\rho_n}\in\SG_n$
from left to right and define the height $h_i$ to be the number of elements $\rho_k$ which are
smaller than $\rho_i$ and lie to the right of $i$, i.e., with $k>i$ and
$\rho_k<\rho_i$.

It is clear that this amounts to another unique encoding of the
permutation $\rho$ by a ``height--vector'' $(h_1,\dots,h_n)$.
This height--vector is translated into a path in the following way.
\begin{quote}
Start at $(0,0)$. 

\noindent
For $i=1,\dots,n$ do the following:
\begin{quote}
If the last point
of the path constructed so far lies below $h_i+1$, then draw as many
up--steps as necessary to reach height $h_i+1$, otherwise make as
may down--jumps as necessary to reach height $h_i+1$.

Finally, draw a
down--step at height $h_i$.
\end{quote}
\end{quote}
For the resulting lattice path $p$, set
$\psi_{(3,1,2)}\of\rho=p$.

\end{dfn}
It is easy to see that $\psi_{(3,1,2)}$ is a well--defined injection.

\begin{dfn}
\label{dfn:psi321}
Define the mapping
$$\psi_{(3,2,1)}:\;\SG_n\rightarrow D^\prime_n$$
by the following construction: Read $\rho=\of{\rho_1,\dots,\rho_n}\in\SG_n$
from left to right and define the height $h_i$ as follows.
\begin{itemize}
\item {\em IF\/} $(i,\rho_i)$ is a \rlmum, then
$h_i$ is the number of elements $\rho_k$ which are smaller than and lie to the right from $\rho_i$, i.e., with $k>i$ and $\rho_k<\rho_i$, 
\item {\em ELSE\/} $h_i$ is the number of elements $\rho_k$ which are bigger than and
lie to the right of $\rho_i$, i.e., with $k>i$ and
$\rho_k>\rho_i$, {\em and\/} which are smaller than the preceding \rlmum.
\end{itemize}

The ``height--vector'' $(h_1,\dots,h_n)$ is
translated into a path in exactly the same way as in
Definition~\ref{dfn:psi312}.
For the resulting lattice path $p$, set
$\psi_{(3,2,1)}\of\rho=p$.

\end{dfn}
It is easy to see that $\psi_{(3,2,1)}$ is a well--defined injection, too.

See Figure~\ref{fig:bijections} for an illustration of both injections
$\psi_{(3,1,2)}$ and $\psi_{(3,2,1)}$ applied to the same permutation
$\of{4,3,5,1,2}$. Note that the entries of the permutation correspond
to the down--steps of the path, as before. Again, we shall more or
less identify a permutation and its corresponding path.

So clearly, our direction to proceed is as follows:
\bit
\item Determine the subsets of $D^\prime$ which
	correspond to permutations containing one, two, etc. occurrences of
	pattern $\of{3,2,1}$ or $\of{3,1,2}$, respectively,
\item identify the generating functions for such subsets.
\eit
In order to describe the relevant subsets of paths, we must have a close
look at the properties of the mappings $\psi_{(3,1,2)}$ and
$\psi_{(3,2,1)}$.

\begin{obs}
\label{obs:simple}
Let $\tau$ be one of the patterns $\of{3,1,2}$ or $\of{3,2,1}$.

\begin{enumerate}

\item
By construction, a jump in $\psi_{\tau}\of{\rho}$ must
be preceded and followed immediately by a down--step; i.e., patterns
\epsffile{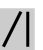}\ and 
\epsffile{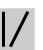}\ are impossible. If a jump occurs
between the $i$--th and $(i+1)$--th down--step, we call it a {\em jump at
position $i$\/}.

\item
For every down--step and for every jump in $\psi_{\tau}\of{\rho}$, there must exist a  preceding \rlmum\ (in the sense of Definition~\ref{dfn:prec-rlmum}).

\item
Every \rlmum\ in permutation $\rho$ corresponds to a ``peak 
\epsffile{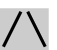}'' in the lattice path $\psi\of\rho$.
However, there might be peaks in $\psi\of\rho$
which do not correspond to \rlma\ in $\rho$ (see Figure~\ref{fig:bijections}
for an example).

\item
Every $\tau$--occurrence $\of{\sigma_1,\sigma_2,\sigma_3}$ in permutaion
$\rho$ is associated with a jump in $\psi\of\rho$: $\sigma_1$ lies weakly
to the left of the jump, $\sigma_2$ and $\sigma_3$ to the right.

\item
$\psi_{\tau}\of\rho$ is an {\em ordinary\/} Dyck--path (i.e.,
	a Dyck--path without any down--jump) if and only if $\rho$ is
	$\tau$--avoiding: Mappings $\psi_{\tau}$ and  $\psi$ coincide, when
	restricted to $\tau$--avoiding permutations, as a second look at the
	proof given in Section~\ref{sec:occ0} shows immediately.

\item
To the right of an arbitrary down--step of height $l$ in
	$\psi_{\tau}\of{\rho}$, there must be at least $l$ down--steps:
	This simple statement implies that to the right of a jump of depth $d$,
	there must be at least $d$ up--steps.
	
\end{enumerate}

\end{obs}

\input fig3

Now we investigate the relationship between jumps in $\psi_\tau\of\rho$
and the number of $\tau$--occurrences in $\rho$ more closely.
\begin{lem}
\label{lem:basic}
Let $\tau\in\SG_3$ be one of the patterns $\of{3,1,2}$ or $\of{3,2,1}$,
let $\rho$ be an arbitrary permutation.
A jump of depth $d$ at position $i$ in $\psi_{\tau}\of{\rho}$,
	which is followed immediately by $l$ ($l>0$ by
	Observation~\ref{obs:simple}!) consecutive down--steps,
	implies $(d\cdot l)$ $\tau$--occurrences of type
	$\of{\rho_r,\rho_{i+j},\rho_{k}}$, where
	\bit
		\item $\rho_r$ is the \rlmum\ preceding the jump, 
		\item $\rho_{i+j}$ are the entries of $\rho$ corresponding to the
			consecutive down--steps immediately following the jump
			($j=1,\dots,l$),
		\item $\rho_k$ is one of the $d$ entries of $\rho$ with $k > i$
			which are 
			\bit
				\item bigger than $\rho_{i+1}$ and smaller than $\rho_i$ for
					$\tau=\of{3,1,2}$,
				\item smaller than $\rho_{i+1}$ for
					$\tau=\of{3,2,1}$.
			\eit
		\eit
\end{lem}

\begin{proof}
The statement sounds a bit complicated, but it is easily
obtained ``by inspection'': 

For the case $\tau=\of{3,1,2}$, consider Figure~\ref{fig:lemma312}
and observe that entry $9$ is the \rlmum\ preceding the jump, that
$\of{2,1}$ are the entries corresponding to the down--steps immediately
following the jump, and that the entries $\of{4,3}$ are the 2 entries to the right of the jump, which are bigger than the entry following the jump and
smaller than the entry preceding the jump --- in some sense, they
``cause'' the jump of depth 2: 
The four occurrences of $\of{3,1,2}$ in Figure~\ref{fig:lemma312}, as
stated in the lemma, are $\of{9,2,4}$, $\of{9,2,3}$,
$\of{9,1,4}$ and $\of{9,1,3}$.

For the case $\tau=\of{3,2,1}$, consider Figure~\ref{fig:lemma321}
and observe that entry $9$ is the \rlmum\ preceding the jump, that
$\of{6,8}$ are the entries corresponding to the down--steps immediately
following the jump, and that entries $\of{4,5}$ are the 2 entries to the
right of the jump, which are smaller than the entry following the jump ---
in the same sense as above, they ``cause'' the jump
of depth 2.
The four occurrences of $\of{3,2,1}$ in Figure~\ref{fig:lemma321}, as
stated in the lemma, are $\of{9,6,4}$, $\of{9,6,5}$,
$\of{9,8,4}$ and $\of{9,8,5}$.
\end{proof}

\input fig4

\input fig5

\begin{cor}
\label{cor:jumpsum}
For $r>0$ and $\tau=\of{3,1,2}$ or $\tau=\of{3,2,1}$, $\psi_{\tau}$ takes
the set $\SGntr{n}{\tau}{r}$ into $\bigcup_{s=1}^r D_{n,s}$; i.e.,
a permutation with exactly $r$ occurrences of $\tau$ is mapped
to a Dyck--path with at most $r$ down--jumps.
\end{cor}

It is immediately obvious from Figures~\ref{fig:lemma312} and
\ref{fig:lemma321} that the enumeration of occurrences associated
with some jump, as given in Lemma~\ref{lem:basic}, is in general incomplete.
Note that down--steps {\em immediately preceding\/}
a jump obviously add to the number of $(3,1,2)$--occurrences, while
this is not the case for $\of{3,2,1}$--occurrences. To the contrary,
we find additional occurrences of pattern $\of{3,2,1}$ if there are
``not enough'' such down--steps.

While the general
situation obviously can get quite complicated, we state the following
extensions to Lemma~\ref{lem:basic}.

\begin{cor}
\label{cor:basic2}
Let $\tau=\of{3,1,2}$, let $\rho$ be an arbitrary permutation.
Consider a jump of depth $d$ at position $i$ in $\psi_{\tau}\of{\rho}$,
which is followed immediately by $l$ consecutive down--steps,
and which is preceded immediately by $m$ consecutive
down--steps {\em including\/} the one corresponding to the preceding
\rlmum. Such a jump implies $(m\cdot d\cdot l)$ $\tau$--occurrences of type
$\of{\rho_{i-g+1},\rho_{i+j},\rho_{k}}$, where
	\bit
		\item $\rho_{i-g+1}$ are the consecutive down--steps immediately
			preceding the jump ($g=1,\dots,m$),
		\item $\rho_{i+j}$ are the entries of $\rho$ corresponding to the
			consecutive down--steps immediately following the jump
			($j=1,\dots,l$),
		\item $\rho_k$ is one of the $d$ entries of $\rho$ with $k > i$
			which are bigger than $\rho_{i+1}$ and smaller than $\rho_i$.
	\eit

\end{cor}

\begin{proof}
In our graphical notation, the situation is as follows:

\bigskip\centerline{\epsffile{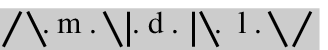}.}\bigskip\noindent
The assertion is easily obtained ``by inspection'', see
Figure~\ref{fig:lemma312}. 
\end{proof}

By careful inspection, we obtain yet another observation:
\begin{cor}
\label{cor:basic3}
In the same situation as in Corollary~\ref{cor:basic2},
consider the \rlma\ $\rho_{i_1},\dots,\rho_{i_r}$
to the left of the \rlmum\ preceding the jump, $\rho_{i-m+1}$. For each
of these, compute the
number $s_j$ of up--steps to the right, which do not correspond to peaks
of \rlma\ and which are to the left of the jump. Let $x$ be the smallest index
for which $s_x < m$. Then there are exactly
\begin{equation*}
d\cdot l\cdot \sum_{\xi=x}^r\of{m-s_\xi}
\end{equation*}
$\tau$--occurrences of type
$\of{\rho_{\xi},\rho_{i+j},\rho_{k}}$, where
	\bit
		\item $\rho_{\xi}$ are the \rlma\ as described above ($\xi=x,\dots,r$),
		\item $\rho_{i+j}$ are the entries of $\rho$ corresponding to the
			consecutive down--steps immediately following the jump
			($j=1,\dots,l$),
		\item $\rho_k$ is one of the $d$ entries of $\rho$ with $k > i$
			which are bigger than $\rho_{i+1}$ and smaller than $\rho_i$.
	\eit

\end{cor}

\begin{lem}
\label{lem:basic2}
Let $\tau=\of{3,2,1}$, let $\rho$ be an arbitrary permutation.
Consider the {\em first\/} jump (counted from the left). Assume that
this jump
\bit
\item occurs at position $i$ in $\psi_{\tau}\of{\rho}$,
\item is of depth $d$,
\item and is followed immediately by $l$ consecutive down--steps.
\eit
Consider the \rlma\ $\rho_{i_1},\dots,\rho_{i_r}$ to the left of the
jump, and their respective heights $h_{i_1},\dots,h_{i_r}$.
For each of these $r$ \rlma\, compute the the number $s_j$ of
down--steps to the right, which do  {\em not correspond\/} to \rlma\ and
which are to the left of the jump.

Let $x$ be the smallest index for which
\begin{equation*}
h_{i_x} - d - s_x > 0.
\end{equation*}
Then there are exactly
\begin{equation*}
d\cdot\sum_{g=x}^{r}\min\of{h_{i_g}- d - s_g, l}
\end{equation*}
$\tau$--occurrences of type $\of{\rho_{i_g},\rho_{i+j},\rho_{k}}$,
where
	\bit
		\item $\rho_{i_g}$ are the \rlma\ as above ($g=1,\dots,r$),
		\item $\rho_{i+j}$ are the entries of $\rho$ corresponding to the
			consecutive down--steps after the jump
			($j=1,\dots,\min\of{h_{i_g} - d - s_g, l}$),
		\item $\rho_k$ are the $d$ entries of $\rho$ with $k > i$
			which are smaller than $\rho_{i+1}$.
	\eit
\end{lem}

\begin{proof}
The statement sounds rather complicated, but again, it is easy to obtain
``by inspection'' --- see Figure~\ref{fig:321complicated}. There, we have
$d=1$ and $l=4$.
The 5 \rlma\ preceding the first (and only)
jump are 6, 8, 10, 11 and
14; with heights 5, 5, 5, 4 and 5; the numbers $s_j$ read 4, 3, 2, 1 and 0.

\input fig6

Observe that the condition for $h_{i_g}$ and $s_g$ simply determines whether
\rlmum\ $\rho_{i_g}$ does ``intersect the jump--configuration'',
which consists of $\of{7,9,12,13,5}$ in Figure~\ref{fig:321complicated}.
Recall that $h_{i_g}$ is the number of entries to the right of
$\rho_{i_g}$ which are
smaller than $\rho_{i_g}$; in order that $\rho_{i_g}$ is involved in a $\tau$--occurrence,
the $d$ entries which ``cause'' the jump must be among them,
{\em and\/} there must be some entry left among them, which is not
``used up'' on the way to the jump. This consideration explains the sum
in the lemma.

In Figure~\ref{fig:321complicated}, the numbers $(h_{i_g}-d-s_g)$ are
$0,1,2,2,4$; corresponding to no $(3,2,1)$--occurrence involving the
first \rlmum\ 6, 1 occurrence $(8,7,5)$, two occurrences $(10,9,5)$
and $(10,7,5)$, another two occurrences $(11,9,5)$ and $(11,7,5)$, and
finally four occurrences $(14,13,5)$, $(14,12,5)$, $(14,9,5)$ and $(14,7,5)$.
\end{proof}

Now we have the prerequisites ready for the remaining proofs. In the following,
we will determine the subsets $P\of{\tau,r}\subseteq D^\prime$, which
correspond bijectively to permutations with exactly $r$ occurrences of $\tau$.
More precisely, we will partition
the set $P\of{\tau,r}$ into subsets $P_{1},\dots,P_{m}$ according to some
``pattern in the paths'', where all paths in $P_{i}$ have the same number
$s_i$ of down--jumps.
Define the weight $w\of p$ of a path in $D_{n,s}$ to be $x^{(2n+s)/2}$
(i.e., assign weight $\sqrt x$ to each up--step and each down--step; and
weight $1$ to each of the $s$ down--jumps), and 
consider the generating function 
$$P_{\tau,r,i}\of x:=\sum_{p}w\of p,$$
where the summation runs over all
$p\in P_{i}$. Then the
generating function we are interested in is given as
\begin{equation}
\label{eq:sliced}
\GFxtr{x}{\tau}{r}=\sum_{i=1}^m x^{-s_i/2} P_{\tau,r,i}\of x.
\end{equation}

%% file: fig1.tex
\begin{figure}
\caption{Illustration of the mapping $\psi$ from permutations to
Dyck--paths. (Note: The permutation $(3,6,5,4,8,2,10,9,7,1)$ is
$(3,1,2)$--avoiding.)}
\label{fig:bijection0}
\begin{picture}(11,11)(0,0)
\put(0.5,10.5){\line(1,0){10}}
\put(0.5,9.5){\line(1,0){10}}
\put(0.5,8.5){\line(1,0){10}}
\put(0.5,7.5){\line(1,0){10}}
\put(0.5,6.5){\line(1,0){10}}
\put(0.5,5.5){\line(1,0){10}}
\put(0.5,4.5){\line(1,0){10}}
\put(0.5,3.5){\line(1,0){10}}
\put(0.5,2.5){\line(1,0){10}}
\put(0.5,1.5){\line(1,0){10}}
\put(0.5,0.5){\line(1,0){10}}
\put(1,0){{\small 1}}
\put(2,0){{\small 2}}
\put(3,0){{\small 3}}
\put(4,0){{\small 4}}
\put(5,0){{\small 5}}
\put(6,0){{\small 6}}
\put(7,0){{\small 7}}
\put(8,0){{\small 8}}
\put(9,0){{\small 9}}
\put(10,0){{\small 10}}
\put(10.5,0.5){\line(0,1){10}}
\put(9.5,0.5){\line(0,1){10}}
\put(8.5,0.5){\line(0,1){10}}
\put(7.5,0.5){\line(0,1){10}}
\put(6.5,0.5){\line(0,1){10}}
\put(5.5,0.5){\line(0,1){10}}
\put(4.5,0.5){\line(0,1){10}}
\put(3.5,0.5){\line(0,1){10}}
\put(2.5,0.5){\line(0,1){10}}
\put(1.5,0.5){\line(0,1){10}}
\put(0.5,0.5){\line(0,1){10}}
\put(0,3){{\small 3}}
\put(0,6){{\small 6}}
\put(0,5){{\small 5}}
\put(0,4){{\small 4}}
\put(0,8){{\small 8}}
\put(0,2){{\small 2}}
\put(0,10){{\small 10}}
\put(0,9){{\small 9}}
\put(0,7){{\small 7}}
\put(0,1){{\small 1}}
\put(1,3){\circle*{0.3}}
\put(2,6){\circle*{0.3}}
\put(3,5){\circle*{0.3}}
\put(4,4){\circle*{0.3}}
\put(5,8){\circle*{0.3}}
\put(6,7){\circle*{0.3}}
\put(7,10){\circle*{0.3}}
\put(8,9){\circle*{0.3}}
\put(9,2){\circle*{0.3}}
\put(10,1){\circle*{0.3}}
\put(0.5,2.5){\usebox{\schraffenbox}}
\put(1.5,2.5){\usebox{\schraffenbox}}
\put(2.5,2.5){\usebox{\schraffenbox}}
\put(3.5,2.5){\usebox{\schraffenbox}}
\put(4.5,2.5){\usebox{\schraffenbox}}
\put(5.5,2.5){\usebox{\schraffenbox}}
\put(6.5,2.5){\usebox{\schraffenbox}}
\put(7.5,2.5){\usebox{\schraffenbox}}
\put(8.5,2.5){\usebox{\schraffenbox}}
\put(9.5,2.5){\usebox{\schraffenbox}}
\put(0.5,0.5){\usebox{\schraffenbox}}
\put(0.5,1.5){\usebox{\schraffenbox}}
\put(0.5,2.5){\usebox{\schraffenbox}}
\put(1.5,5.5){\usebox{\schraffenbox}}
\put(2.5,5.5){\usebox{\schraffenbox}}
\put(3.5,5.5){\usebox{\schraffenbox}}
\put(4.5,5.5){\usebox{\schraffenbox}}
\put(5.5,5.5){\usebox{\schraffenbox}}
\put(6.5,5.5){\usebox{\schraffenbox}}
\put(7.5,5.5){\usebox{\schraffenbox}}
\put(8.5,5.5){\usebox{\schraffenbox}}
\put(9.5,5.5){\usebox{\schraffenbox}}
\put(1.5,0.5){\usebox{\schraffenbox}}
\put(1.5,1.5){\usebox{\schraffenbox}}
\put(1.5,2.5){\usebox{\schraffenbox}}
\put(1.5,3.5){\usebox{\schraffenbox}}
\put(1.5,4.5){\usebox{\schraffenbox}}
\put(1.5,5.5){\usebox{\schraffenbox}}
\put(4.5,7.5){\usebox{\schraffenbox}}
\put(5.5,7.5){\usebox{\schraffenbox}}
\put(6.5,7.5){\usebox{\schraffenbox}}
\put(7.5,7.5){\usebox{\schraffenbox}}
\put(8.5,7.5){\usebox{\schraffenbox}}
\put(9.5,7.5){\usebox{\schraffenbox}}
\put(4.5,0.5){\usebox{\schraffenbox}}
\put(4.5,1.5){\usebox{\schraffenbox}}
\put(4.5,2.5){\usebox{\schraffenbox}}
\put(4.5,3.5){\usebox{\schraffenbox}}
\put(4.5,4.5){\usebox{\schraffenbox}}
\put(4.5,5.5){\usebox{\schraffenbox}}
\put(4.5,6.5){\usebox{\schraffenbox}}
\put(4.5,7.5){\usebox{\schraffenbox}}
\put(6.5,9.5){\usebox{\schraffenbox}}
\put(7.5,9.5){\usebox{\schraffenbox}}
\put(8.5,9.5){\usebox{\schraffenbox}}
\put(9.5,9.5){\usebox{\schraffenbox}}
\put(6.5,0.5){\usebox{\schraffenbox}}
\put(6.5,1.5){\usebox{\schraffenbox}}
\put(6.5,2.5){\usebox{\schraffenbox}}
\put(6.5,3.5){\usebox{\schraffenbox}}
\put(6.5,4.5){\usebox{\schraffenbox}}
\put(6.5,5.5){\usebox{\schraffenbox}}
\put(6.5,6.5){\usebox{\schraffenbox}}
\put(6.5,7.5){\usebox{\schraffenbox}}
\put(6.5,8.5){\usebox{\schraffenbox}}
\put(6.5,9.5){\usebox{\schraffenbox}}
\linethickness{2pt}
\put(1.5,3.5){\line(0,1){3}}
\put(1.5,6.5){\line(1,0){1}}
\put(2,6){\circle{0.6}}
\put(2.5,6.5){\line(1,0){1}}
\put(3.5,6.5){\line(1,0){1}}
\put(4.5,6.5){\line(0,1){2}}
\put(4.5,8.5){\line(1,0){1}}
\put(5,8){\circle{0.6}}
\put(0.5,0.5){\line(0,1){3}}
\put(0.5,3.5){\line(1,0){1}}
\put(1,3){\circle{0.6}}
\put(7.5,10.5){\line(1,0){1}}
\put(8.5,10.5){\line(1,0){1}}
\put(9.5,10.5){\line(1,0){1}}
\put(5.5,8.5){\line(1,0){1}}
\put(6.5,8.5){\line(0,1){2}}
\put(6.5,10.5){\line(1,0){1}}
\put(7,10){\circle{0.6}}
\end{picture}

\setlength{\unitlength}{0.6cm}
\begin{picture}(21,6)(-0.5,-0.5)
\put(0,0){\line(1,0){20}}
\put(0,1){\line(1,0){20}}
\put(0,2){\line(1,0){20}}
\put(0,3){\line(1,0){20}}
\put(0,4){\line(1,0){20}}
\put(0,5){\line(1,0){20}}
\put(0,0){\line(0,1){5}}
\put(1,0){\line(0,1){5}}
\put(2,0){\line(0,1){5}}
\put(3,0){\line(0,1){5}}
\put(4,0){\line(0,1){5}}
\put(5,0){\line(0,1){5}}
\put(6,0){\line(0,1){5}}
\put(7,0){\line(0,1){5}}
\put(8,0){\line(0,1){5}}
\put(9,0){\line(0,1){5}}
\put(10,0){\line(0,1){5}}
\put(11,0){\line(0,1){5}}
\put(12,0){\line(0,1){5}}
\put(13,0){\line(0,1){5}}
\put(14,0){\line(0,1){5}}
\put(15,0){\line(0,1){5}}
\put(16,0){\line(0,1){5}}
\put(17,0){\line(0,1){5}}
\put(18,0){\line(0,1){5}}
\put(19,0){\line(0,1){5}}
\put(20,0){\line(0,1){5}}
\put(0,0){\line(1,1){3}}
\put(3,3){\line(1,-1){1}}
\put(3.6,2.6){{\small 3}}
\put(4,2){\line(1,1){3}}
\put(7,5){\line(1,-1){1}}
\put(7.6,4.6){{\small 6}}
\put(8,4){\line(1,-1){1}}
\put(8.6,3.6){{\small 5}}
\put(9,3){\line(1,-1){1}}
\put(9.6,2.6){{\small 4}}
\put(10,2){\line(1,1){2}}
\put(12,4){\line(1,-1){1}}
\put(12.6,3.6){{\small 8}}
\put(13,3){\line(1,-1){1}}
\put(13.6,2.6){{\small 7}}
\put(14,2){\line(1,1){2}}
\put(16,4){\line(1,-1){1}}
\put(16.6,3.6){{\small 10}}
\put(17,3){\line(1,-1){1}}
\put(17.6,2.6){{\small 9}}
\put(18,2){\line(1,-1){1}}
\put(18.6,1.6){{\small 2}}
\put(19,1){\line(1,-1){1}}
\put(19.6,0.6){{\small 1}}
\end{picture}

\end{figure}
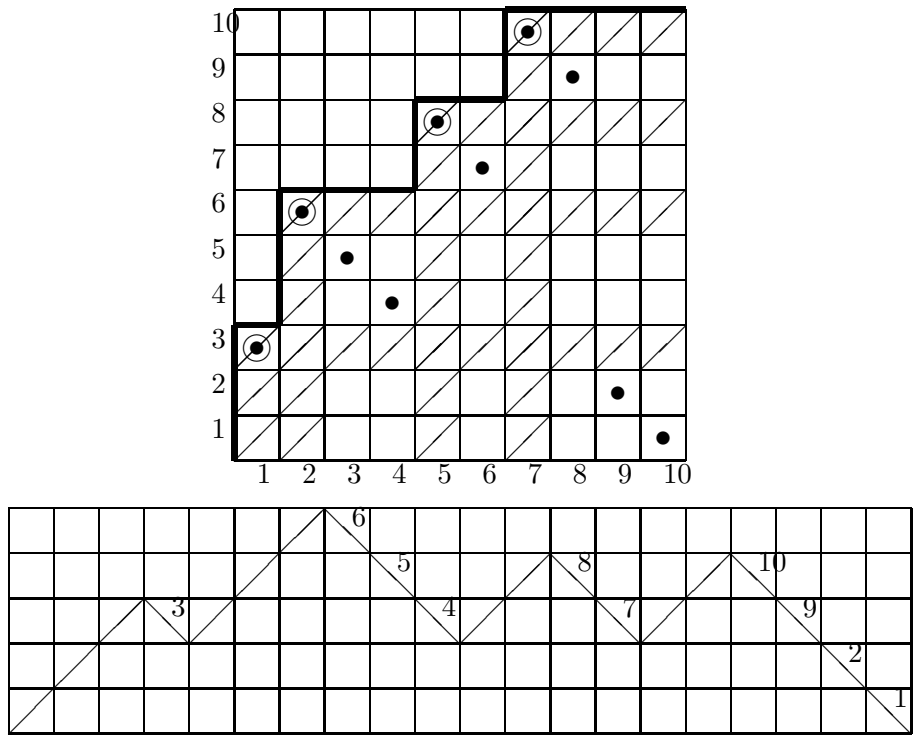

%% file: fig2.tex
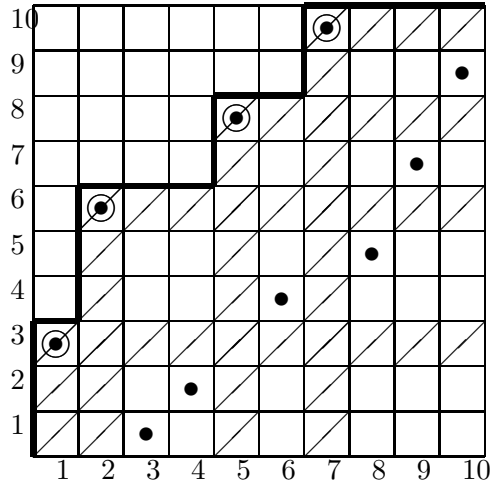
\begin{figure}
\caption{Illustration of the mapping $\psi$ from permutations to
Dyck--paths. (Note: The permutation $(3,6,1,2,8,4,10,5,7,9)$ is
$(3,2,1)$--avoiding.)}
\label{fig:bijection0-1}
\begin{picture}(11,11)(0,0)
\put(0.5,10.5){\line(1,0){10}}
\put(0.5,9.5){\line(1,0){10}}
\put(0.5,8.5){\line(1,0){10}}
\put(0.5,7.5){\line(1,0){10}}
\put(0.5,6.5){\line(1,0){10}}
\put(0.5,5.5){\line(1,0){10}}
\put(0.5,4.5){\line(1,0){10}}
\put(0.5,3.5){\line(1,0){10}}
\put(0.5,2.5){\line(1,0){10}}
\put(0.5,1.5){\line(1,0){10}}
\put(0.5,0.5){\line(1,0){10}}
\put(1,0){{\small 1}}
\put(2,0){{\small 2}}
\put(3,0){{\small 3}}
\put(4,0){{\small 4}}
\put(5,0){{\small 5}}
\put(6,0){{\small 6}}
\put(7,0){{\small 7}}
\put(8,0){{\small 8}}
\put(9,0){{\small 9}}
\put(10,0){{\small 10}}
\put(10.5,0.5){\line(0,1){10}}
\put(9.5,0.5){\line(0,1){10}}
\put(8.5,0.5){\line(0,1){10}}
\put(7.5,0.5){\line(0,1){10}}
\put(6.5,0.5){\line(0,1){10}}
\put(5.5,0.5){\line(0,1){10}}
\put(4.5,0.5){\line(0,1){10}}
\put(3.5,0.5){\line(0,1){10}}
\put(2.5,0.5){\line(0,1){10}}
\put(1.5,0.5){\line(0,1){10}}
\put(0.5,0.5){\line(0,1){10}}
\put(0,3){{\small 3}}
\put(0,6){{\small 6}}
\put(0,1){{\small 1}}
\put(0,2){{\small 2}}
\put(0,8){{\small 8}}
\put(0,4){{\small 4}}
\put(0,10){{\small 10}}
\put(0,5){{\small 5}}
\put(0,7){{\small 7}}
\put(0,9){{\small 9}}
\put(1,3){\circle*{0.3}}
\put(2,6){\circle*{0.3}}
\put(3,1){\circle*{0.3}}
\put(4,2){\circle*{0.3}}
\put(5,8){\circle*{0.3}}
\put(6,4){\circle*{0.3}}
\put(7,10){\circle*{0.3}}
\put(8,5){\circle*{0.3}}
\put(9,7){\circle*{0.3}}
\put(10,9){\circle*{0.3}}
\put(0.5,0.5){\line(0,1){3}}
\put(0.5,3.5){\line(1,0){1}}
\put(1,3){\circle{0.6}}
\put(0.5,2.5){\usebox{\schraffenbox}}
\put(1.5,2.5){\usebox{\schraffenbox}}
\put(2.5,2.5){\usebox{\schraffenbox}}
\put(3.5,2.5){\usebox{\schraffenbox}}
\put(4.5,2.5){\usebox{\schraffenbox}}
\put(5.5,2.5){\usebox{\schraffenbox}}
\put(6.5,2.5){\usebox{\schraffenbox}}
\put(7.5,2.5){\usebox{\schraffenbox}}
\put(8.5,2.5){\usebox{\schraffenbox}}
\put(9.5,2.5){\usebox{\schraffenbox}}
\put(0.5,0.5){\usebox{\schraffenbox}}
\put(0.5,1.5){\usebox{\schraffenbox}}
\put(0.5,2.5){\usebox{\schraffenbox}}
\put(1.5,3.5){\line(0,1){3}}
\put(1.5,6.5){\line(1,0){1}}
\put(2,6){\circle{0.6}}
\put(1.5,5.5){\usebox{\schraffenbox}}
\put(2.5,5.5){\usebox{\schraffenbox}}
\put(3.5,5.5){\usebox{\schraffenbox}}
\put(4.5,5.5){\usebox{\schraffenbox}}
\put(5.5,5.5){\usebox{\schraffenbox}}
\put(6.5,5.5){\usebox{\schraffenbox}}
\put(7.5,5.5){\usebox{\schraffenbox}}
\put(8.5,5.5){\usebox{\schraffenbox}}
\put(9.5,5.5){\usebox{\schraffenbox}}
\put(1.5,0.5){\usebox{\schraffenbox}}
\put(1.5,1.5){\usebox{\schraffenbox}}
\put(1.5,2.5){\usebox{\schraffenbox}}
\put(1.5,3.5){\usebox{\schraffenbox}}
\put(1.5,4.5){\usebox{\schraffenbox}}
\put(1.5,5.5){\usebox{\schraffenbox}}
\put(2.5,6.5){\line(1,0){1}}
\put(3.5,6.5){\line(1,0){1}}
\put(4.5,6.5){\line(0,1){2}}
\put(4.5,8.5){\line(1,0){1}}
\put(5,8){\circle{0.6}}
\put(4.5,7.5){\usebox{\schraffenbox}}
\put(5.5,7.5){\usebox{\schraffenbox}}
\put(6.5,7.5){\usebox{\schraffenbox}}
\put(7.5,7.5){\usebox{\schraffenbox}}
\put(8.5,7.5){\usebox{\schraffenbox}}
\put(9.5,7.5){\usebox{\schraffenbox}}
\put(4.5,0.5){\usebox{\schraffenbox}}
\put(4.5,1.5){\usebox{\schraffenbox}}
\put(4.5,2.5){\usebox{\schraffenbox}}
\put(4.5,3.5){\usebox{\schraffenbox}}
\put(4.5,4.5){\usebox{\schraffenbox}}
\put(4.5,5.5){\usebox{\schraffenbox}}
\put(4.5,6.5){\usebox{\schraffenbox}}
\put(4.5,7.5){\usebox{\schraffenbox}}
\put(5.5,8.5){\line(1,0){1}}
\put(6.5,8.5){\line(0,1){2}}
\put(6.5,10.5){\line(1,0){1}}
\put(7,10){\circle{0.6}}
\put(6.5,9.5){\usebox{\schraffenbox}}
\put(7.5,9.5){\usebox{\schraffenbox}}
\put(8.5,9.5){\usebox{\schraffenbox}}
\put(9.5,9.5){\usebox{\schraffenbox}}
\put(6.5,0.5){\usebox{\schraffenbox}}
\put(6.5,1.5){\usebox{\schraffenbox}}
\put(6.5,2.5){\usebox{\schraffenbox}}
\put(6.5,3.5){\usebox{\schraffenbox}}
\put(6.5,4.5){\usebox{\schraffenbox}}
\put(6.5,5.5){\usebox{\schraffenbox}}
\put(6.5,6.5){\usebox{\schraffenbox}}
\put(6.5,7.5){\usebox{\schraffenbox}}
\put(6.5,8.5){\usebox{\schraffenbox}}
\put(6.5,9.5){\usebox{\schraffenbox}}
\linethickness{2pt}
\put(1.5,3.5){\line(0,1){3}}
\put(1.5,6.5){\line(1,0){1}}
\put(2,6){\circle{0.6}}
\put(2.5,6.5){\line(1,0){1}}
\put(3.5,6.5){\line(1,0){1}}
\put(4.5,6.5){\line(0,1){2}}
\put(4.5,8.5){\line(1,0){1}}
\put(5,8){\circle{0.6}}
\put(0.5,0.5){\line(0,1){3}}
\put(0.5,3.5){\line(1,0){1}}
\put(1,3){\circle{0.6}}
\put(7.5,10.5){\line(1,0){1}}
\put(8.5,10.5){\line(1,0){1}}
\put(9.5,10.5){\line(1,0){1}}
\put(5.5,8.5){\line(1,0){1}}
\put(6.5,8.5){\line(0,1){2}}
\put(6.5,10.5){\line(1,0){1}}
\put(7,10){\circle{0.6}}
\end{picture}
\end{figure}

%% file: fig3.tex
\begin{figure}
\caption{Illustration of the bijections $\psi_{(3,2,1)}$ and
$\psi_{(3,2,1)}$ for the permutation $\rho=(4,3,5,1,2)$.}
\label{fig:bijections}
\begin{picture}(6,6)(0,0)
\put(0.5,5.5){\line(1,0){5}}
\put(0.5,4.5){\line(1,0){5}}
\put(0.5,3.5){\line(1,0){5}}
\put(0.5,2.5){\line(1,0){5}}
\put(0.5,1.5){\line(1,0){5}}
\put(0.5,0.5){\line(1,0){5}}
\put(1,0){{\small 1}}
\put(2,0){{\small 2}}
\put(3,0){{\small 3}}
\put(4,0){{\small 4}}
\put(5,0){{\small 5}}
\put(5.5,0.5){\line(0,1){5}}
\put(4.5,0.5){\line(0,1){5}}
\put(3.5,0.5){\line(0,1){5}}
\put(2.5,0.5){\line(0,1){5}}
\put(1.5,0.5){\line(0,1){5}}
\put(0.5,0.5){\line(0,1){5}}
\put(0,4){{\small 4}}
\put(0,3){{\small 3}}
\put(0,5){{\small 5}}
\put(0,1){{\small 1}}
\put(0,2){{\small 2}}
\put(1,4){\circle*{0.3}}
\put(2,3){\circle*{0.3}}
\put(3,5){\circle*{0.3}}
\put(4,1){\circle*{0.3}}
\put(5,2){\circle*{0.3}}
\put(1,4){\circle{0.6}}
\put(3,5){\circle{0.6}}
\put(0.5,3.5){\usebox{\schraffenbox}}
\put(1.5,3.5){\usebox{\schraffenbox}}
\put(2.5,3.5){\usebox{\schraffenbox}}
\put(3.5,3.5){\usebox{\schraffenbox}}
\put(4.5,3.5){\usebox{\schraffenbox}}
\put(0.5,0.5){\usebox{\schraffenbox}}
\put(0.5,1.5){\usebox{\schraffenbox}}
\put(0.5,2.5){\usebox{\schraffenbox}}
\put(0.5,3.5){\usebox{\schraffenbox}}
\put(2.5,4.5){\usebox{\schraffenbox}}
\put(3.5,4.5){\usebox{\schraffenbox}}
\put(4.5,4.5){\usebox{\schraffenbox}}
\put(2.5,0.5){\usebox{\schraffenbox}}
\put(2.5,1.5){\usebox{\schraffenbox}}
\put(2.5,2.5){\usebox{\schraffenbox}}
\put(2.5,3.5){\usebox{\schraffenbox}}
\put(2.5,4.5){\usebox{\schraffenbox}}
\thicklines
\put(0.5,0.5){\line(0,1){4}}
\put(0.5,4.5){\line(1,0){1}}
\put(1.5,4.5){\line(1,0){1}}
\put(2.5,4.5){\line(0,1){1}}
\put(2.5,5.5){\line(1,0){1}}
\put(3.5,5.5){\line(1,0){1}}
\put(4.5,5.5){\line(1,0){1}}
\end{picture}

\begin{center}
The image of $\rho$ under $\psi_{(3,1,2)}$:  One jump of depth 1.
\end{center}

\setlength{\unitlength}{0.6cm}

\begin{picture}(12,5)(-0.5,-0.5)
\put(0,0){\line(1,0){11}}
\put(0,1){\line(1,0){11}}
\put(0,2){\line(1,0){11}}
\put(0,3){\line(1,0){11}}
\put(0,4){\line(1,0){11}}
\put(0,0){\line(0,1){4}}
\put(1,0){\line(0,1){4}}
\put(2,0){\line(0,1){4}}
\put(3,0){\line(0,1){4}}
\put(4,0){\line(0,1){4}}
\put(5,0){\line(0,1){4}}
\put(6,0){\line(0,1){4}}
\put(7,0){\line(0,1){4}}
\put(8,0){\line(0,1){4}}
\put(9,0){\line(0,1){4}}
\put(10,0){\line(0,1){4}}
\put(11,0){\line(0,1){4}}
\put(-0.4,-0.150000){{\small 0}}
\put(-0.4,0.850000){{\small 1}}
\put(-0.4,1.850000){{\small 2}}
\put(-0.4,2.850000){{\small 3}}
\put(-0.4,3.850000){{\small 4}}
\put(-0.100000,-0.5){{\small 0}}
\put(0.900000,-0.5){{\small 1}}
\put(1.900000,-0.5){{\small 2}}
\put(2.900000,-0.5){{\small 3}}
\put(3.900000,-0.5){{\small 4}}
\put(4.900000,-0.5){{\small 5}}
\put(5.900000,-0.5){{\small 6}}
\put(6.900000,-0.5){{\small 7}}
\put(7.900000,-0.5){{\small 8}}
\put(8.900000,-0.5){{\small 9}}
\put(9.900000,-0.5){{\small 10}}
\put(10.900000,-0.5){{\small 11}}
\put(0,0){\line(1,1){4}}
\put(4,4){\line(1,-1){1}}
\put(5,3){\line(1,-1){1}}
\put(6,2){\line(1,1){1}}
\put(7,3){\line(1,-1){1}}
\put(8,1){\line(1,-1){1}}
\put(9,0){\line(1,1){1}}
\put(10,1){\line(1,-1){1}}
\put(4.6,3.6){{\small 4}}
\put(5.6,2.6){{\small 3}}
\put(7.6,2.6){{\small 5}}
\put(8.6,0.6){{\small 1}}
\put(10.6,0.6){{\small 2}}
\end{picture}

\begin{center}
The image of $\rho$ under $\psi_{(3,2,1)}$: One jump of depth 2.
\end{center}

\begin{picture}(13,5)(-0.5,-0.5)
\put(0,0){\line(1,0){12}}
\put(0,1){\line(1,0){12}}
\put(0,2){\line(1,0){12}}
\put(0,3){\line(1,0){12}}
\put(0,4){\line(1,0){12}}
\put(0,0){\line(0,1){4}}
\put(1,0){\line(0,1){4}}
\put(2,0){\line(0,1){4}}
\put(3,0){\line(0,1){4}}
\put(4,0){\line(0,1){4}}
\put(5,0){\line(0,1){4}}
\put(6,0){\line(0,1){4}}
\put(7,0){\line(0,1){4}}
\put(8,0){\line(0,1){4}}
\put(9,0){\line(0,1){4}}
\put(10,0){\line(0,1){4}}
\put(11,0){\line(0,1){4}}
\put(12,0){\line(0,1){4}}
\put(-0.4,-0.150000){{\small 0}}
\put(-0.4,0.850000){{\small 1}}
\put(-0.4,1.850000){{\small 2}}
\put(-0.4,2.850000){{\small 3}}
\put(-0.4,3.850000){{\small 4}}
\put(-0.100000,-0.5){{\small 0}}
\put(0.900000,-0.5){{\small 1}}
\put(1.900000,-0.5){{\small 2}}
\put(2.900000,-0.5){{\small 3}}
\put(3.900000,-0.5){{\small 4}}
\put(4.900000,-0.5){{\small 5}}
\put(5.900000,-0.5){{\small 6}}
\put(6.900000,-0.5){{\small 7}}
\put(7.900000,-0.5){{\small 8}}
\put(8.900000,-0.5){{\small 9}}
\put(9.900000,-0.5){{\small 10}}
\put(10.900000,-0.5){{\small 11}}
\put(11.900000,-0.5){{\small 12}}
\put(0,0){\line(1,1){4}}
\put(4,4){\line(1,-1){1}}
\put(5,1){\line(1,-1){1}}
\put(6,0){\line(1,1){3}}
\put(9,3){\line(1,-1){1}}
\put(10,2){\line(1,-1){1}}
\put(11,1){\line(1,-1){1}}
\put(4.6,3.6){{\small 4}}
\put(5.6,0.6){{\small 3}}
\put(9.6,2.6){{\small 5}}
\put(10.6,1.6){{\small 1}}
\put(11.6,0.6){{\small 2}}
\end{picture}
\end{figure}
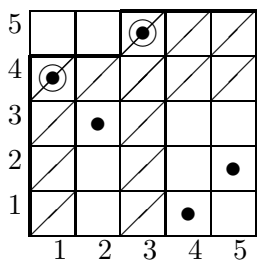

%% file: fig4.tex
\begin{figure}
\caption{Illustration for Lemma~\ref{lem:basic}. (Note: The permutation $(7,9,8,6,5,2,1,4,3)$ is not
$(3,1,2)$--avoiding.)}
\label{fig:lemma312}
\begin{picture}(10,10)(0,0)
\put(0.5,9.5){\line(1,0){9}}
\put(0.5,8.5){\line(1,0){9}}
\put(0.5,7.5){\line(1,0){9}}
\put(0.5,6.5){\line(1,0){9}}
\put(0.5,5.5){\line(1,0){9}}
\put(0.5,4.5){\line(1,0){9}}
\put(0.5,3.5){\line(1,0){9}}
\put(0.5,2.5){\line(1,0){9}}
\put(0.5,1.5){\line(1,0){9}}
\put(0.5,0.5){\line(1,0){9}}
\put(1,0){{\small 1}}
\put(2,0){{\small 2}}
\put(3,0){{\small 3}}
\put(4,0){{\small 4}}
\put(5,0){{\small 5}}
\put(6,0){{\small 6}}
\put(7,0){{\small 7}}
\put(8,0){{\small 8}}
\put(9,0){{\small 9}}
\put(9.5,0.5){\line(0,1){9}}
\put(8.5,0.5){\line(0,1){9}}
\put(7.5,0.5){\line(0,1){9}}
\put(6.5,0.5){\line(0,1){9}}
\put(5.5,0.5){\line(0,1){9}}
\put(4.5,0.5){\line(0,1){9}}
\put(3.5,0.5){\line(0,1){9}}
\put(2.5,0.5){\line(0,1){9}}
\put(1.5,0.5){\line(0,1){9}}
\put(0.5,0.5){\line(0,1){9}}
\put(0,7){{\small 7}}
\put(0,9){{\small 9}}
\put(0,8){{\small 8}}
\put(0,6){{\small 6}}
\put(0,5){{\small 5}}
\put(0,2){{\small 2}}
\put(0,1){{\small 1}}
\put(0,4){{\small 4}}
\put(0,3){{\small 3}}
\put(1,7){\circle*{0.3}}
\put(2,9){\circle*{0.3}}
\put(3,8){\circle*{0.3}}
\put(4,6){\circle*{0.3}}
\put(5,5){\circle*{0.3}}
\put(6,2){\circle*{0.3}}
\put(7,1){\circle*{0.3}}
\put(8,4){\circle*{0.3}}
\put(9,3){\circle*{0.3}}
\put(1,7){\circle{0.6}}
\put(2,9){\circle{0.6}}
\put(0.5,6.5){\usebox{\schraffenbox}}
\put(1.5,6.5){\usebox{\schraffenbox}}
\put(2.5,6.5){\usebox{\schraffenbox}}
\put(3.5,6.5){\usebox{\schraffenbox}}
\put(4.5,6.5){\usebox{\schraffenbox}}
\put(5.5,6.5){\usebox{\schraffenbox}}
\put(6.5,6.5){\usebox{\schraffenbox}}
\put(7.5,6.5){\usebox{\schraffenbox}}
\put(8.5,6.5){\usebox{\schraffenbox}}
\put(0.5,0.5){\usebox{\schraffenbox}}
\put(0.5,1.5){\usebox{\schraffenbox}}
\put(0.5,2.5){\usebox{\schraffenbox}}
\put(0.5,3.5){\usebox{\schraffenbox}}
\put(0.5,4.5){\usebox{\schraffenbox}}
\put(0.5,5.5){\usebox{\schraffenbox}}
\put(0.5,6.5){\usebox{\schraffenbox}}
\put(1.5,8.5){\usebox{\schraffenbox}}
\put(2.5,8.5){\usebox{\schraffenbox}}
\put(3.5,8.5){\usebox{\schraffenbox}}
\put(4.5,8.5){\usebox{\schraffenbox}}
\put(5.5,8.5){\usebox{\schraffenbox}}
\put(6.5,8.5){\usebox{\schraffenbox}}
\put(7.5,8.5){\usebox{\schraffenbox}}
\put(8.5,8.5){\usebox{\schraffenbox}}
\put(1.5,0.5){\usebox{\schraffenbox}}
\put(1.5,1.5){\usebox{\schraffenbox}}
\put(1.5,2.5){\usebox{\schraffenbox}}
\put(1.5,3.5){\usebox{\schraffenbox}}
\put(1.5,4.5){\usebox{\schraffenbox}}
\put(1.5,5.5){\usebox{\schraffenbox}}
\put(1.5,6.5){\usebox{\schraffenbox}}
\put(1.5,7.5){\usebox{\schraffenbox}}
\put(1.5,8.5){\usebox{\schraffenbox}}
\thicklines
\put(0.5,0.5){\line(0,1){7}}
\put(0.5,7.5){\line(1,0){1}}
\put(1.5,7.5){\line(0,1){2}}
\put(1.5,9.5){\line(1,0){1}}
\put(2.5,9.5){\line(1,0){1}}
\put(3.5,9.5){\line(1,0){1}}
\put(4.5,9.5){\line(1,0){1}}
\put(5.5,9.5){\line(1,0){1}}
\put(6.5,9.5){\line(1,0){1}}
\put(7.5,9.5){\line(1,0){1}}
\put(8.5,9.5){\line(1,0){1}}
\end{picture}

\setlength{\unitlength}{0.6cm}
\begin{picture}(21,9)(-0.5,-0.5)
\put(0,0){\line(1,0){20}}
\put(0,1){\line(1,0){20}}
\put(0,2){\line(1,0){20}}
\put(0,3){\line(1,0){20}}
\put(0,4){\line(1,0){20}}
\put(0,5){\line(1,0){20}}
\put(0,6){\line(1,0){20}}
\put(0,7){\line(1,0){20}}
\put(0,8){\line(1,0){20}}
\put(0,0){\line(0,1){8}}
\put(1,0){\line(0,1){8}}
\put(2,0){\line(0,1){8}}
\put(3,0){\line(0,1){8}}
\put(4,0){\line(0,1){8}}
\put(5,0){\line(0,1){8}}
\put(6,0){\line(0,1){8}}
\put(7,0){\line(0,1){8}}
\put(8,0){\line(0,1){8}}
\put(9,0){\line(0,1){8}}
\put(10,0){\line(0,1){8}}
\put(11,0){\line(0,1){8}}
\put(12,0){\line(0,1){8}}
\put(13,0){\line(0,1){8}}
\put(14,0){\line(0,1){8}}
\put(15,0){\line(0,1){8}}
\put(16,0){\line(0,1){8}}
\put(17,0){\line(0,1){8}}
\put(18,0){\line(0,1){8}}
\put(19,0){\line(0,1){8}}
\put(20,0){\line(0,1){8}}
\put(-0.4,-0.150000){{\small 0}}
\put(-0.4,0.850000){{\small 1}}
\put(-0.4,1.850000){{\small 2}}
\put(-0.4,2.850000){{\small 3}}
\put(-0.4,3.850000){{\small 4}}
\put(-0.4,4.850000){{\small 5}}
\put(-0.4,5.850000){{\small 6}}
\put(-0.4,6.850000){{\small 7}}
\put(-0.4,7.850000){{\small 8}}
\put(-0.100000,-0.5){{\small 0}}
\put(0.900000,-0.5){{\small 1}}
\put(1.900000,-0.5){{\small 2}}
\put(2.900000,-0.5){{\small 3}}
\put(3.900000,-0.5){{\small 4}}
\put(4.900000,-0.5){{\small 5}}
\put(5.900000,-0.5){{\small 6}}
\put(6.900000,-0.5){{\small 7}}
\put(7.900000,-0.5){{\small 8}}
\put(8.900000,-0.5){{\small 9}}
\put(9.900000,-0.5){{\small 10}}
\put(10.900000,-0.5){{\small 11}}
\put(11.900000,-0.5){{\small 12}}
\put(12.900000,-0.5){{\small 13}}
\put(13.900000,-0.5){{\small 14}}
\put(14.900000,-0.5){{\small 15}}
\put(15.900000,-0.5){{\small 16}}
\put(16.900000,-0.5){{\small 17}}
\put(17.900000,-0.5){{\small 18}}
\put(18.900000,-0.5){{\small 19}}
\put(19.900000,-0.5){{\small 20}}
\put(0,0){\line(1,1){7}}
\put(7,7){\line(1,-1){1}}
\put(8,6){\line(1,1){2}}
\put(10,8){\line(1,-1){1}}
\put(11,7){\line(1,-1){1}}
\put(12,6){\line(1,-1){1}}
\put(13,5){\line(1,-1){1}}
\put(14,2){\line(1,-1){1}}
\put(15,1){\line(1,-1){1}}
\put(16,0){\line(1,1){2}}
\put(18,2){\line(1,-1){1}}
\put(19,1){\line(1,-1){1}}
\put(7.6,6.6){{\small 7}}
\put(10.6,7.6){{\small 9}}
\put(11.6,6.6){{\small 8}}
\put(12.6,5.6){{\small 6}}
\put(13.6,4.6){{\small 5}}
\put(14.6,1.6){{\small 2}}
\put(15.6,0.6){{\small 1}}
\put(18.6,1.6){{\small 4}}
\put(19.6,0.6){{\small 3}}
\end{picture}

\end{figure}

%% file: fig5.tex
\begin{figure}
\caption{Illustration for Lemma~\ref{lem:basic}. (Note: The permutation $(7,1,2,9,3,6,8,4,5)$ is not
$(3,2,1)$--avoiding.)}
\label{fig:lemma321}
\begin{picture}(10,10)(0,0)
\put(0.5,9.5){\line(1,0){9}}
\put(0.5,8.5){\line(1,0){9}}
\put(0.5,7.5){\line(1,0){9}}
\put(0.5,6.5){\line(1,0){9}}
\put(0.5,5.5){\line(1,0){9}}
\put(0.5,4.5){\line(1,0){9}}
\put(0.5,3.5){\line(1,0){9}}
\put(0.5,2.5){\line(1,0){9}}
\put(0.5,1.5){\line(1,0){9}}
\put(0.5,0.5){\line(1,0){9}}
\put(1,0){{\small 1}}
\put(2,0){{\small 2}}
\put(3,0){{\small 3}}
\put(4,0){{\small 4}}
\put(5,0){{\small 5}}
\put(6,0){{\small 6}}
\put(7,0){{\small 7}}
\put(8,0){{\small 8}}
\put(9,0){{\small 9}}
\put(9.5,0.5){\line(0,1){9}}
\put(8.5,0.5){\line(0,1){9}}
\put(7.5,0.5){\line(0,1){9}}
\put(6.5,0.5){\line(0,1){9}}
\put(5.5,0.5){\line(0,1){9}}
\put(4.5,0.5){\line(0,1){9}}
\put(3.5,0.5){\line(0,1){9}}
\put(2.5,0.5){\line(0,1){9}}
\put(1.5,0.5){\line(0,1){9}}
\put(0.5,0.5){\line(0,1){9}}
\put(0,7){{\small 7}}
\put(0,1){{\small 1}}
\put(0,2){{\small 2}}
\put(0,9){{\small 9}}
\put(0,3){{\small 3}}
\put(0,6){{\small 6}}
\put(0,8){{\small 8}}
\put(0,4){{\small 4}}
\put(0,5){{\small 5}}
\put(1,7){\circle*{0.3}}
\put(2,1){\circle*{0.3}}
\put(3,2){\circle*{0.3}}
\put(4,9){\circle*{0.3}}
\put(5,3){\circle*{0.3}}
\put(6,6){\circle*{0.3}}
\put(7,8){\circle*{0.3}}
\put(8,4){\circle*{0.3}}
\put(9,5){\circle*{0.3}}
\put(1,7){\circle{0.6}}
\put(4,9){\circle{0.6}}
\put(0.5,6.5){\usebox{\schraffenbox}}
\put(1.5,6.5){\usebox{\schraffenbox}}
\put(2.5,6.5){\usebox{\schraffenbox}}
\put(3.5,6.5){\usebox{\schraffenbox}}
\put(4.5,6.5){\usebox{\schraffenbox}}
\put(5.5,6.5){\usebox{\schraffenbox}}
\put(6.5,6.5){\usebox{\schraffenbox}}
\put(7.5,6.5){\usebox{\schraffenbox}}
\put(8.5,6.5){\usebox{\schraffenbox}}
\put(0.5,0.5){\usebox{\schraffenbox}}
\put(0.5,1.5){\usebox{\schraffenbox}}
\put(0.5,2.5){\usebox{\schraffenbox}}
\put(0.5,3.5){\usebox{\schraffenbox}}
\put(0.5,4.5){\usebox{\schraffenbox}}
\put(0.5,5.5){\usebox{\schraffenbox}}
\put(0.5,6.5){\usebox{\schraffenbox}}
\put(3.5,8.5){\usebox{\schraffenbox}}
\put(4.5,8.5){\usebox{\schraffenbox}}
\put(5.5,8.5){\usebox{\schraffenbox}}
\put(6.5,8.5){\usebox{\schraffenbox}}
\put(7.5,8.5){\usebox{\schraffenbox}}
\put(8.5,8.5){\usebox{\schraffenbox}}
\put(3.5,0.5){\usebox{\schraffenbox}}
\put(3.5,1.5){\usebox{\schraffenbox}}
\put(3.5,2.5){\usebox{\schraffenbox}}
\put(3.5,3.5){\usebox{\schraffenbox}}
\put(3.5,4.5){\usebox{\schraffenbox}}
\put(3.5,5.5){\usebox{\schraffenbox}}
\put(3.5,6.5){\usebox{\schraffenbox}}
\put(3.5,7.5){\usebox{\schraffenbox}}
\put(3.5,8.5){\usebox{\schraffenbox}}
\thicklines
\put(0.5,0.5){\line(0,1){7}}
\put(0.5,7.5){\line(1,0){1}}
\put(1.5,7.5){\line(1,0){1}}
\put(2.5,7.5){\line(1,0){1}}
\put(3.5,7.5){\line(0,1){2}}
\put(3.5,9.5){\line(1,0){1}}
\put(4.5,9.5){\line(1,0){1}}
\put(5.5,9.5){\line(1,0){1}}
\put(6.5,9.5){\line(1,0){1}}
\put(7.5,9.5){\line(1,0){1}}
\put(8.5,9.5){\line(1,0){1}}
\end{picture}

\setlength{\unitlength}{0.6cm}

\begin{picture}(21,8)(-0.5,-0.5)
\put(0,0){\line(1,0){20}}
\put(0,1){\line(1,0){20}}
\put(0,2){\line(1,0){20}}
\put(0,3){\line(1,0){20}}
\put(0,4){\line(1,0){20}}
\put(0,5){\line(1,0){20}}
\put(0,6){\line(1,0){20}}
\put(0,7){\line(1,0){20}}
\put(0,0){\line(0,1){7}}
\put(1,0){\line(0,1){7}}
\put(2,0){\line(0,1){7}}
\put(3,0){\line(0,1){7}}
\put(4,0){\line(0,1){7}}
\put(5,0){\line(0,1){7}}
\put(6,0){\line(0,1){7}}
\put(7,0){\line(0,1){7}}
\put(8,0){\line(0,1){7}}
\put(9,0){\line(0,1){7}}
\put(10,0){\line(0,1){7}}
\put(11,0){\line(0,1){7}}
\put(12,0){\line(0,1){7}}
\put(13,0){\line(0,1){7}}
\put(14,0){\line(0,1){7}}
\put(15,0){\line(0,1){7}}
\put(16,0){\line(0,1){7}}
\put(17,0){\line(0,1){7}}
\put(18,0){\line(0,1){7}}
\put(19,0){\line(0,1){7}}
\put(20,0){\line(0,1){7}}
\put(-0.4,-0.150000){{\small 0}}
\put(-0.4,0.850000){{\small 1}}
\put(-0.4,1.850000){{\small 2}}
\put(-0.4,2.850000){{\small 3}}
\put(-0.4,3.850000){{\small 4}}
\put(-0.4,4.850000){{\small 5}}
\put(-0.4,5.850000){{\small 6}}
\put(-0.4,6.850000){{\small 7}}
\put(-0.100000,-0.5){{\small 0}}
\put(0.900000,-0.5){{\small 1}}
\put(1.900000,-0.5){{\small 2}}
\put(2.900000,-0.5){{\small 3}}
\put(3.900000,-0.5){{\small 4}}
\put(4.900000,-0.5){{\small 5}}
\put(5.900000,-0.5){{\small 6}}
\put(6.900000,-0.5){{\small 7}}
\put(7.900000,-0.5){{\small 8}}
\put(8.900000,-0.5){{\small 9}}
\put(9.900000,-0.5){{\small 10}}
\put(10.900000,-0.5){{\small 11}}
\put(11.900000,-0.5){{\small 12}}
\put(12.900000,-0.5){{\small 13}}
\put(13.900000,-0.5){{\small 14}}
\put(14.900000,-0.5){{\small 15}}
\put(15.900000,-0.5){{\small 16}}
\put(16.900000,-0.5){{\small 17}}
\put(17.900000,-0.5){{\small 18}}
\put(18.900000,-0.5){{\small 19}}
\put(19.900000,-0.5){{\small 20}}
\put(0,0){\line(1,1){7}}
\put(7,7){\line(1,-1){1}}
\put(8,6){\line(1,-1){1}}
\put(9,5){\line(1,-1){1}}
\put(10,4){\line(1,1){2}}
\put(12,6){\line(1,-1){1}}
\put(13,5){\line(1,-1){1}}
\put(14,2){\line(1,-1){1}}
\put(15,1){\line(1,-1){1}}
\put(16,0){\line(1,1){2}}
\put(18,2){\line(1,-1){1}}
\put(19,1){\line(1,-1){1}}
\put(7.6,6.6){{\small 7}}
\put(8.6,5.6){{\small 1}}
\put(9.6,4.6){{\small 2}}
\put(12.6,5.6){{\small 9}}
\put(13.6,4.6){{\small 3}}
\put(14.6,1.6){{\small 6}}
\put(15.6,0.6){{\small 8}}
\put(18.6,1.6){{\small 4}}
\put(19.6,0.6){{\small 5}}
\end{picture}

\end{figure}

%% file: fig6.tex
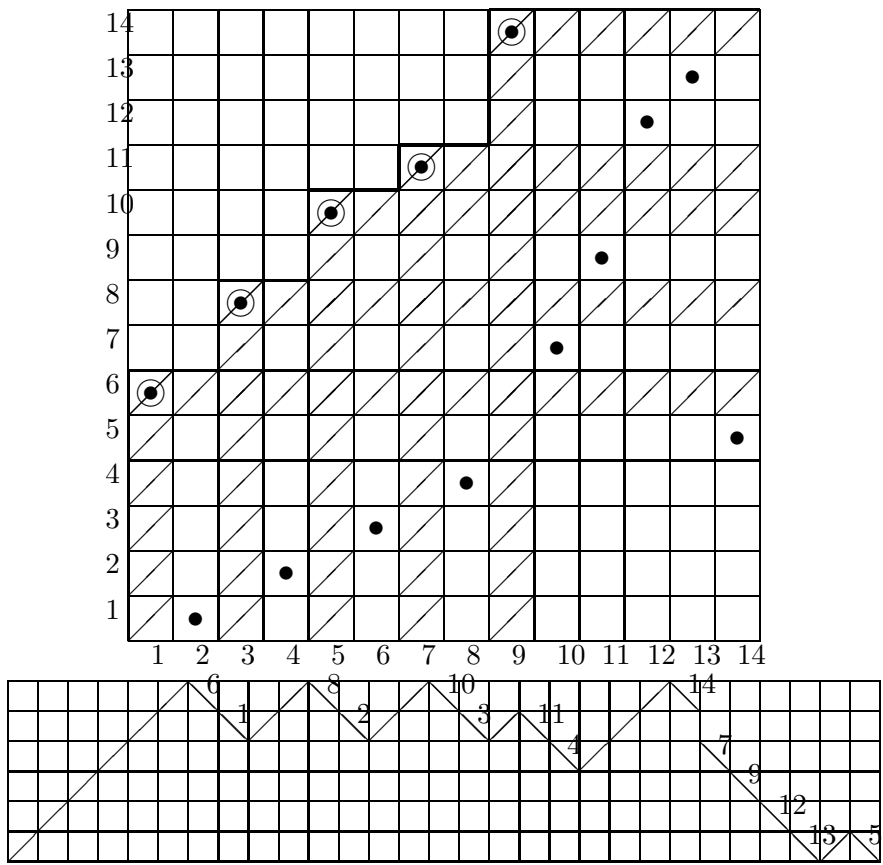
\begin{figure}
\caption{Illustration for Lemma~\ref{lem:basic2}.}
\label{fig:321complicated}
\begin{picture}(15,15)(0,0)
\put(0.5,14.5){\line(1,0){14}}
\put(0.5,13.5){\line(1,0){14}}
\put(0.5,12.5){\line(1,0){14}}
\put(0.5,11.5){\line(1,0){14}}
\put(0.5,10.5){\line(1,0){14}}
\put(0.5,9.5){\line(1,0){14}}
\put(0.5,8.5){\line(1,0){14}}
\put(0.5,7.5){\line(1,0){14}}
\put(0.5,6.5){\line(1,0){14}}
\put(0.5,5.5){\line(1,0){14}}
\put(0.5,4.5){\line(1,0){14}}
\put(0.5,3.5){\line(1,0){14}}
\put(0.5,2.5){\line(1,0){14}}
\put(0.5,1.5){\line(1,0){14}}
\put(0.5,0.5){\line(1,0){14}}
\put(1,0){{\small 1}}
\put(2,0){{\small 2}}
\put(3,0){{\small 3}}
\put(4,0){{\small 4}}
\put(5,0){{\small 5}}
\put(6,0){{\small 6}}
\put(7,0){{\small 7}}
\put(8,0){{\small 8}}
\put(9,0){{\small 9}}
\put(10,0){{\small 10}}
\put(11,0){{\small 11}}
\put(12,0){{\small 12}}
\put(13,0){{\small 13}}
\put(14,0){{\small 14}}
\put(14.5,0.5){\line(0,1){14}}
\put(13.5,0.5){\line(0,1){14}}
\put(12.5,0.5){\line(0,1){14}}
\put(11.5,0.5){\line(0,1){14}}
\put(10.5,0.5){\line(0,1){14}}
\put(9.5,0.5){\line(0,1){14}}
\put(8.5,0.5){\line(0,1){14}}
\put(7.5,0.5){\line(0,1){14}}
\put(6.5,0.5){\line(0,1){14}}
\put(5.5,0.5){\line(0,1){14}}
\put(4.5,0.5){\line(0,1){14}}
\put(3.5,0.5){\line(0,1){14}}
\put(2.5,0.5){\line(0,1){14}}
\put(1.5,0.5){\line(0,1){14}}
\put(0.5,0.5){\line(0,1){14}}
\put(0,6){{\small 6}}
\put(0,1){{\small 1}}
\put(0,8){{\small 8}}
\put(0,2){{\small 2}}
\put(0,10){{\small 10}}
\put(0,3){{\small 3}}
\put(0,11){{\small 11}}
\put(0,4){{\small 4}}
\put(0,14){{\small 14}}
\put(0,7){{\small 7}}
\put(0,9){{\small 9}}
\put(0,12){{\small 12}}
\put(0,13){{\small 13}}
\put(0,5){{\small 5}}
\put(1,6){\circle*{0.3}}
\put(2,1){\circle*{0.3}}
\put(3,8){\circle*{0.3}}
\put(4,2){\circle*{0.3}}
\put(5,10){\circle*{0.3}}
\put(6,3){\circle*{0.3}}
\put(7,11){\circle*{0.3}}
\put(8,4){\circle*{0.3}}
\put(9,14){\circle*{0.3}}
\put(10,7){\circle*{0.3}}
\put(11,9){\circle*{0.3}}
\put(12,12){\circle*{0.3}}
\put(13,13){\circle*{0.3}}
\put(14,5){\circle*{0.3}}
\put(1,6){\circle{0.6}}
\put(3,8){\circle{0.6}}
\put(5,10){\circle{0.6}}
\put(7,11){\circle{0.6}}
\put(9,14){\circle{0.6}}
\put(0.5,5.5){\usebox{\schraffenbox}}
\put(1.5,5.5){\usebox{\schraffenbox}}
\put(2.5,5.5){\usebox{\schraffenbox}}
\put(3.5,5.5){\usebox{\schraffenbox}}
\put(4.5,5.5){\usebox{\schraffenbox}}
\put(5.5,5.5){\usebox{\schraffenbox}}
\put(6.5,5.5){\usebox{\schraffenbox}}
\put(7.5,5.5){\usebox{\schraffenbox}}
\put(8.5,5.5){\usebox{\schraffenbox}}
\put(9.5,5.5){\usebox{\schraffenbox}}
\put(10.5,5.5){\usebox{\schraffenbox}}
\put(11.5,5.5){\usebox{\schraffenbox}}
\put(12.5,5.5){\usebox{\schraffenbox}}
\put(13.5,5.5){\usebox{\schraffenbox}}
\put(0.5,0.5){\usebox{\schraffenbox}}
\put(0.5,1.5){\usebox{\schraffenbox}}
\put(0.5,2.5){\usebox{\schraffenbox}}
\put(0.5,3.5){\usebox{\schraffenbox}}
\put(0.5,4.5){\usebox{\schraffenbox}}
\put(0.5,5.5){\usebox{\schraffenbox}}
\put(2.5,7.5){\usebox{\schraffenbox}}
\put(3.5,7.5){\usebox{\schraffenbox}}
\put(4.5,7.5){\usebox{\schraffenbox}}
\put(5.5,7.5){\usebox{\schraffenbox}}
\put(6.5,7.5){\usebox{\schraffenbox}}
\put(7.5,7.5){\usebox{\schraffenbox}}
\put(8.5,7.5){\usebox{\schraffenbox}}
\put(9.5,7.5){\usebox{\schraffenbox}}
\put(10.5,7.5){\usebox{\schraffenbox}}
\put(11.5,7.5){\usebox{\schraffenbox}}
\put(12.5,7.5){\usebox{\schraffenbox}}
\put(13.5,7.5){\usebox{\schraffenbox}}
\put(2.5,0.5){\usebox{\schraffenbox}}
\put(2.5,1.5){\usebox{\schraffenbox}}
\put(2.5,2.5){\usebox{\schraffenbox}}
\put(2.5,3.5){\usebox{\schraffenbox}}
\put(2.5,4.5){\usebox{\schraffenbox}}
\put(2.5,5.5){\usebox{\schraffenbox}}
\put(2.5,6.5){\usebox{\schraffenbox}}
\put(2.5,7.5){\usebox{\schraffenbox}}
\put(4.5,9.5){\usebox{\schraffenbox}}
\put(5.5,9.5){\usebox{\schraffenbox}}
\put(6.5,9.5){\usebox{\schraffenbox}}
\put(7.5,9.5){\usebox{\schraffenbox}}
\put(8.5,9.5){\usebox{\schraffenbox}}
\put(9.5,9.5){\usebox{\schraffenbox}}
\put(10.5,9.5){\usebox{\schraffenbox}}
\put(11.5,9.5){\usebox{\schraffenbox}}
\put(12.5,9.5){\usebox{\schraffenbox}}
\put(13.5,9.5){\usebox{\schraffenbox}}
\put(4.5,0.5){\usebox{\schraffenbox}}
\put(4.5,1.5){\usebox{\schraffenbox}}
\put(4.5,2.5){\usebox{\schraffenbox}}
\put(4.5,3.5){\usebox{\schraffenbox}}
\put(4.5,4.5){\usebox{\schraffenbox}}
\put(4.5,5.5){\usebox{\schraffenbox}}
\put(4.5,6.5){\usebox{\schraffenbox}}
\put(4.5,7.5){\usebox{\schraffenbox}}
\put(4.5,8.5){\usebox{\schraffenbox}}
\put(4.5,9.5){\usebox{\schraffenbox}}
\put(6.5,10.5){\usebox{\schraffenbox}}
\put(7.5,10.5){\usebox{\schraffenbox}}
\put(8.5,10.5){\usebox{\schraffenbox}}
\put(9.5,10.5){\usebox{\schraffenbox}}
\put(10.5,10.5){\usebox{\schraffenbox}}
\put(11.5,10.5){\usebox{\schraffenbox}}
\put(12.5,10.5){\usebox{\schraffenbox}}
\put(13.5,10.5){\usebox{\schraffenbox}}
\put(6.5,0.5){\usebox{\schraffenbox}}
\put(6.5,1.5){\usebox{\schraffenbox}}
\put(6.5,2.5){\usebox{\schraffenbox}}
\put(6.5,3.5){\usebox{\schraffenbox}}
\put(6.5,4.5){\usebox{\schraffenbox}}
\put(6.5,5.5){\usebox{\schraffenbox}}
\put(6.5,6.5){\usebox{\schraffenbox}}
\put(6.5,7.5){\usebox{\schraffenbox}}
\put(6.5,8.5){\usebox{\schraffenbox}}
\put(6.5,9.5){\usebox{\schraffenbox}}
\put(6.5,10.5){\usebox{\schraffenbox}}
\put(8.5,13.5){\usebox{\schraffenbox}}
\put(9.5,13.5){\usebox{\schraffenbox}}
\put(10.5,13.5){\usebox{\schraffenbox}}
\put(11.5,13.5){\usebox{\schraffenbox}}
\put(12.5,13.5){\usebox{\schraffenbox}}
\put(13.5,13.5){\usebox{\schraffenbox}}
\put(8.5,0.5){\usebox{\schraffenbox}}
\put(8.5,1.5){\usebox{\schraffenbox}}
\put(8.5,2.5){\usebox{\schraffenbox}}
\put(8.5,3.5){\usebox{\schraffenbox}}
\put(8.5,4.5){\usebox{\schraffenbox}}
\put(8.5,5.5){\usebox{\schraffenbox}}
\put(8.5,6.5){\usebox{\schraffenbox}}
\put(8.5,7.5){\usebox{\schraffenbox}}
\put(8.5,8.5){\usebox{\schraffenbox}}
\put(8.5,9.5){\usebox{\schraffenbox}}
\put(8.5,10.5){\usebox{\schraffenbox}}
\put(8.5,11.5){\usebox{\schraffenbox}}
\put(8.5,12.5){\usebox{\schraffenbox}}
\put(8.5,13.5){\usebox{\schraffenbox}}
\thicklines
\put(0.5,0.5){\line(0,1){6}}
\put(0.5,6.5){\line(1,0){1}}
\put(1.5,6.5){\line(1,0){1}}
\put(2.5,6.5){\line(0,1){2}}
\put(2.5,8.5){\line(1,0){1}}
\put(3.5,8.5){\line(1,0){1}}
\put(4.5,8.5){\line(0,1){2}}
\put(4.5,10.5){\line(1,0){1}}
\put(5.5,10.5){\line(1,0){1}}
\put(6.5,10.5){\line(0,1){1}}
\put(6.5,11.5){\line(1,0){1}}
\put(7.5,11.5){\line(1,0){1}}
\put(8.5,11.5){\line(0,1){3}}
\put(8.5,14.5){\line(1,0){1}}
\put(9.5,14.5){\line(1,0){1}}
\put(10.5,14.5){\line(1,0){1}}
\put(11.5,14.5){\line(1,0){1}}
\put(12.5,14.5){\line(1,0){1}}
\put(13.5,14.5){\line(1,0){1}}
\end{picture}

\setlength{\unitlength}{0.4cm}
\begin{picture}(30,7)(-0.5,-0.5)
\put(0,0){\line(1,0){29}}
\put(0,1){\line(1,0){29}}
\put(0,2){\line(1,0){29}}
\put(0,3){\line(1,0){29}}
\put(0,4){\line(1,0){29}}
\put(0,5){\line(1,0){29}}
\put(0,6){\line(1,0){29}}
\put(0,0){\line(0,1){6}}
\put(1,0){\line(0,1){6}}
\put(2,0){\line(0,1){6}}
\put(3,0){\line(0,1){6}}
\put(4,0){\line(0,1){6}}
\put(5,0){\line(0,1){6}}
\put(6,0){\line(0,1){6}}
\put(7,0){\line(0,1){6}}
\put(8,0){\line(0,1){6}}
\put(9,0){\line(0,1){6}}
\put(10,0){\line(0,1){6}}
\put(11,0){\line(0,1){6}}
\put(12,0){\line(0,1){6}}
\put(13,0){\line(0,1){6}}
\put(14,0){\line(0,1){6}}
\put(15,0){\line(0,1){6}}
\put(16,0){\line(0,1){6}}
\put(17,0){\line(0,1){6}}
\put(18,0){\line(0,1){6}}
\put(19,0){\line(0,1){6}}
\put(20,0){\line(0,1){6}}
\put(21,0){\line(0,1){6}}
\put(22,0){\line(0,1){6}}
\put(23,0){\line(0,1){6}}
\put(24,0){\line(0,1){6}}
\put(25,0){\line(0,1){6}}
\put(26,0){\line(0,1){6}}
\put(27,0){\line(0,1){6}}
\put(28,0){\line(0,1){6}}
\put(29,0){\line(0,1){6}}
\put(0,0){\line(1,1){6}}
\put(6,6){\line(1,-1){1}}
\put(7,5){\line(1,-1){1}}
\put(8,4){\line(1,1){2}}
\put(10,6){\line(1,-1){1}}
\put(11,5){\line(1,-1){1}}
\put(12,4){\line(1,1){2}}
\put(14,6){\line(1,-1){1}}
\put(15,5){\line(1,-1){1}}
\put(16,4){\line(1,1){1}}
\put(17,5){\line(1,-1){1}}
\put(18,4){\line(1,-1){1}}
\put(19,3){\line(1,1){3}}
\put(22,6){\line(1,-1){1}}
\put(23,4){\line(1,-1){1}}
\put(24,3){\line(1,-1){1}}
\put(25,2){\line(1,-1){1}}
\put(26,1){\line(1,-1){1}}
\put(27,0){\line(1,1){1}}
\put(28,1){\line(1,-1){1}}
\put(6.6,5.6){{\small 6}}
\put(7.6,4.6){{\small 1}}
\put(10.6,5.6){{\small 8}}
\put(11.6,4.6){{\small 2}}
\put(14.6,5.6){{\small 10}}
\put(15.6,4.6){{\small 3}}
\put(17.6,4.6){{\small 11}}
\put(18.6,3.6){{\small 4}}
\put(22.6,5.6){{\small 14}}
\put(23.6,3.6){{\small 7}}
\put(24.6,2.6){{\small 9}}
\put(25.6,1.6){{\small 12}}
\put(26.6,0.6){{\small 13}}
\put(28.6,0.6){{\small 5}}
\end{picture}

\end{figure}

%% file: proof1.tex
\subsection{The case of one occurrence of length--$3$--patterns}
\label{sec:occ1}
What we want to do is ``assemble the respective lattice paths
from appropriate parts'', and ``translate'' this construction into the
corresponding generating function. For convenience, we introduce the notation
$$ c = \frac{1}{2x}-\frac{1}{2x}\sqrt{1-4x} $$
for the generating function of ordinary Dyck--paths (without jumps).

Observe that lattice paths from $(a,0)$ to
$(b,l)$ ($b\geq a$, $l \geq 0$), which do not go below the horizontal axis,
can be ``composed'' quite easily (see Figure~\ref{fig:genfunc}),
and that the respective generating function is 
\begin{equation}
\label{eq:path0-l}
c^{l+1}x^{l/2}. 
\end{equation}
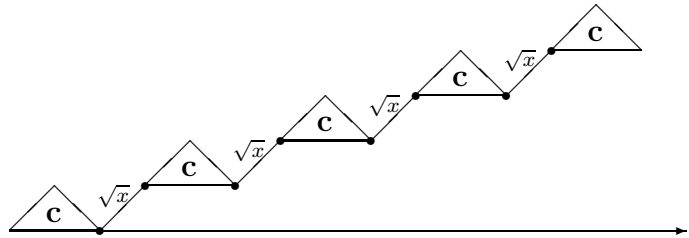
\begin{figure}
\caption{The generating function for lattice paths ``climbing up'' from
height 0 to height 4, which never go below the horizontal axis, is
$c^5 \left(\sqrt{x}\right)^4=c^5 x^2$.}
\label{fig:genfunc}

\begin{picture}(13,7)(-1,0)
\put(-1,1){\vector(1,0){15}}
\put(-1,1){\usebox{\Dyckbox}}
\put(1,1){\usebox{\upstepbox}}
\put(2,2){\usebox{\Dyckbox}}
\put(4,2){\usebox{\upstepbox}}
\put(5,3){\usebox{\Dyckbox}}
\put(7,3){\usebox{\upstepbox}}
\put(8,4){\usebox{\Dyckbox}}
\put(10,4){\usebox{\upstepbox}}
\put(11,5){\usebox{\Dyckbox}}
\end{picture}
\end{figure}

Now let us turn to our lattice paths with jumps:
By Corollary~\ref{cor:jumpsum}, for an arbitrary permutation $\rho$ with
exactly one occurrence of $\tau=\of{3,1,2}$ or $\tau=\of{3,2,1}$, respectively,
the paths $\psi_\tau\of{\rho}$ must contain exactly one
down--jump.  Moreover, by Corollaries~\ref{cor:basic2} and \ref{cor:basic3} or Lemma~\ref{lem:basic2},
respectively, there must hold certain conditions on
\bit
\item down--steps {\em before\/},
\item down--steps {\em after\/},
\item and \rlma\ {\em before\/}
\eit
this single jump.

Let us start with a proof of \eqref{eq:S312-1}: Observe that for any
permutation $\rho$ with exactly one occurrence of $\tau=\of{3,1,2}$, the
path--segment around the jump must look {\em exactly\/} like
\epsffile{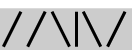}\ 
by Corollaries~\ref{cor:basic2} and \ref{cor:basic3}.
In turn, it is easy to see that
all lattice paths with precisely one such ``down--jump path--segment''
and no other jumps
correspond to permutations with precisely one occurrence
of $\of{3,1,2}$.

So we can immediately write down the generating function $\GFxtr{x}{(3,1,2)}{1}$ we are interested in:
\begin{equation*}
\GFxtr{x}{(3,1,2)}{1} =
	\frac{1}{\sqrt x}
	\sum_{l=1}^{\infty}\of{c^{l+1}{x^{l/2}}}^2x^{5/2}
	= \frac{c^4x^3}{1-c^2x},
\end{equation*}
which coincides with \eqref{eq:S312-1}. (Here, the summation index $l$ is
simply
to be interpreted as the height 
of the path--segment 
\epsffile{p5.eps}.)\qed

\bigskip

Now we turn to the proof of \eqref{eq:S321-1}. The conditions of
Lemma~\ref{lem:basic2} imply, that in order {\em not\/} to have another
\rlmum\ ``intersecting'' our jump (and thus introducing at least a second
occurrence of $\tau=\of{3,2,1}$), there must be at least $\max(1,l)$
down--steps after the \rlmum\ preceding the jump, where $l$ is the height
of the ``valley'' preceding this \rlmum\ --- expressed in our graphical
notation, it is easy to see that we must have precisely the following
situation
\epsffile{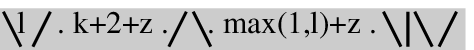}\  
(where $z\geq0$ is an arbitrary integer), at the end of which we are at height $(k+1)$. Consider the generating function
of the ``lattice paths, starting from this path--segment up to the endpoint'':
For $l>0$, it is given as
\begin{equation*}
\sum_{l=1}^{\infty} \frac{x^{(k+5+l)/2}}{1-x} c^{k+2} x^{(k+1)/2} =
\frac{c^2x^{(l+6)/2}}{(1-c x)(1-x)}.
\end{equation*}
(The summands are not simplified on purpose, in order to make ``visible'' the connection
to the ``path--segment''.) For $l=0$, we have
\begin{equation*}
\frac{c^2x^{7/2}}{(1-c x)(1-x)}.
\end{equation*}
Thus, we may write down immediately the desired generating function
$\GFxtr{x}{(3,2,1)}{1}$:
\begin{multline*}
\GFxtr{x}{(3,2,1)}{1} = \\
	\frac{1}{\sqrt x}
	\of{c\frac{c^2x^{7/2}}{(1-c x)(1-x)}+
	\sum_{l=1}^{\infty}
		c^{l+2}x^{(l+1)/2}\frac{c^2x^{(l+6)/2}}{(1-c x)(1-x)}
	} = \\
	\frac{c^3x^3\of{c^2 x-c x + 1}}{(1-x)(1-c x)^2},
\end{multline*}
which coincides with \eqref{eq:S321-1}.\qed

%% file: proof2.tex
\subsection{The case of two occurrences of length--$3$--patterns}
\label{sec:occ2}

\input taubase312

\input taubase321

If a permutation contains exactly two occurrences of $\tau$, then the
corresponding lattice path $\psi_\rho\of{\tau}$ has one or two down--jumps
(by Corollary~\ref{cor:jumpsum}). So we have to consider three cases for
the ``list of depths of jumps'', $\jv$:
\bit
\item $\jv = (1)$, i.e., there is precisely one jump of depth 1,
	(which is ``responsible'' for two $\tau$--occurrences),
\item $\jv = (2)$, i.e., there is precisely one jump of depth 2,
	(which is ``responsible'' for two $\tau$--occurrences),
\item $\jv = (1,1)$, i.e., there are precisely two jumps of depth 1
	each (each of which is ``responsible'' for precisely one
	$\tau$--occurrence).
\eit

As we shall see, the cases must be divided even further: However,
we can make use of certain symmetries.

\begin{dfn}
For a permutation $\rho\in\SG_n$, consider all entries that are involved
in a $\tau$--occurrence. Assume that the length of the subword $\sigma$
formed by these entries is $k$, then there is a permutation $\mu\in\SG_k$ with
entries ``in the same relative order'' as the the entries of $\sigma$: Call
this permutation $\mu$ the {\em \tbase\ of $\rho$\/}.
\end{dfn}
Clearly, the \tbase\ of $\rho$ has the same number of occurrences
as $\rho$ itself. For permutations with exactly one $\tau$--occurrence,
the only possible \tbase\ is $\tau$ itself.

\begin{obs}
For all $r$, the \tbase s of permutations with exactly $r$
$\tau$--occurrences are of length $\leq 3r$: In particular, for $r$
fixed there is only a finite set of possible \tbase s with
$r$ $\tau$--occurrences.

It is immediately clear that the number of $\tau$--occurrences in some
permutation $\rho$ stays unchanged, if
\bit
\item The permutation graph $G\of\rho$ is reflected at the
	downwards--sloping diagonal for both $\tau=\of{3,1,2}$ and
	$\tau=\of{3,2,1}$,
\item The permututation graph $G\of\rho$ is reflected at the
	the upwards--sloping diagonal for $\tau=\of{3,2,1}$ only.
\eit
Therefore, the generating function of permutations with some
fixed \tbase\ is the same as the generating function of permutations
with the ``reflected'' \tbase. 
\end{obs}

Tables~\ref{tab:312} and \ref{tab:321} list all $\tau$--bases for
1 and 2 occurrences of $(3,1,2)$ and $(3,2,1)$, respectively. The
lists were compiled by a brute--force computer search.

\subsubsection{$\tau=\of{3,1,2}$, $r=2$, $\jv=(2)$}
By the same type of considerations as in our proof of
\eqref{eq:S312-1}, we observe that permutations with exactly one jump of depth 2 and exactly two occurrences of $\tau=\of{3,1,2}$ are in one--to--one
correspondence to lattice paths, where the
``path--segment'' around the jump looks {\em exactly\/} like
\epsffile{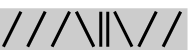}. (The \tbase\ corresponding to these
permutations is $\of{4,1,3,2}$; see Table~\ref{tab:312} e.)
We can immediately write down the generating function of such lattice
paths:
\begin{equation}
\label{eq:S312-2-2}
\sum_{l=1}^{\infty}\of{c^{l+1}{x^{l/2}}}x^{7/2}\of{c^{l+2}{x^{(l+1)/2}}}
	= \frac{c^5x^5}{1-c^2x}.
\end{equation}
Here, the summation index $l$ should be be interpreted as the height where
the path--segment begins. The summands are not simplified on purpose, in order
to make visible the ``composition'' of the paths.

\subsubsection{$\tau=\of{3,1,2}$, $r=2$, $\jv=(1)$} Recall the
considerations for the proof of \eqref{eq:S312-1}: 
There are exactly three ways to obtain precisely two occurrences from a
single jump of depth 1, namely
\begin{enumerate}
\item There is a second ``peak'' immediately preceding the jump, i.e.,
	we face path--segment 
	\epsffile{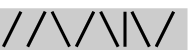}\ (\tbase\ $\of{3,4,1,2}$, see
		Table~\ref{tab:312} d),
\item There is a second down--step immediately before the jump, i.e.,
	we face path--segment 
	\epsffile{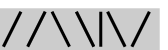}\ (\tbase\ $\of{4,3,1,2}$, see
		Table~\ref{tab:312} h),
\item There is a second down--step immediately after the jump, i.e.,
	we face path--segment 
	\epsffile{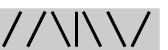}\  (\tbase\ $\of{4,2,1,3}$, see
		Table~\ref{tab:312} f).
\end{enumerate}
The generating function for the permutations corresponding to the
first case is $x \GFxtr{x}{\of{3,1,2}}{1}$: We simply have to ``insert''
another peak 
\epsffile{ud.eps}\ immediately before the peak preceding the
jump.

For the second case, note that the path--segment 
\epsffile{p9.eps}\ 
can be easily transformed into the path--segment 
\epsffile{p8.eps}\ of the
first case (by replacing a ``peak'' by a ``valley''), hence
the generating function is $x \GFxtr{x}{\of{3,1,2}}{1}$ again.

For the third case, note that the \tbase\ $\of{4,2,1,3}$, when reflected at the downwards--sloping diagonal,
gives $\of{4,1,3,2}$, which corresponds to one jump of depth 2: So
the generating function must be the same as in \eqref{eq:S312-2-2}.

\subsubsection{$\tau=\of{3,1,2}$, $r=2$, $\jv=(1,1)$}
We need a simple generalization of \eqref{eq:path0-l}: The generating
function of all lattice paths which
\bit
\item start at $(a,k)$ for some $a\in{\mathbb Z}$ and $k\geq 0$,
\item end at $(b, l)$ for $l\geq k$ and $b\geq a$,
\item and do not go below the horizontal axis
\eit
can be expressed in terms of $c$ and $x$:
\begin{equation}
\label{eq:pathk-l}
c_{k,l}:= \sum_{h=0}^k x^{(l-k+2h)/2} c^{l-k+2h+1} =
	\frac{\left((c^2x)^{k+1}-1\right)c^{l-k+1}x^{(l-k)/2}}{c^2x-1}.
\end{equation}
(To see this, interpret the summation index $h$ as the {\em minimal\/} height
reached by the respective path.)

Now it is easy to write down the generating function for the set of lattice
paths which do contain precisely two path--segments, connected by some
arbitrary path--segment $p$:

\bigskip\centerline{\epsffile{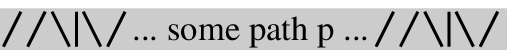}}\bigskip\noindent
(the corresponding \tbase s:
$\of{3,1,2,6,4,5}$,
$\of{3,1,5,2,4}$,
$\of{3,1,6,4,5,2}$ and
$\of{4,2,3,6,1,5}$; see Table~\ref{tab:312} a, b, c and g).
Using \eqref{eq:pathk-l}, we obtain:
\begin{multline}
\label{eq:S312-2-11}
\sum_{l=1}^{\infty}c^{l+1}x^{l/2}x^{5/2}\of{
	\sum_{k=1}^{l}c_{k,l}x^{5/2}c^{k+1}x^{k/2}
	+\sum_{k=l+1}^{\infty}c_{l,k}x^{5/2}c^{k+1}x^{k/2}
	} =\\
\frac{c^5x^6\of{1+c^2x-c^4x^2}}{\of{1-c^2x}^3}.
\end{multline}

\subsubsection{All together}
So we may deduce immediately
\begin{multline*}
\GFxtr{x}{\of{3,1,2}}{2} = 2x \frac{c^4 x^3}{1-c^2x} +
	\frac{2c^5x^5}{x\of{1-c^2x}} +
	\frac{c^5x^6\of{1+c^2x-c^4x^2}}{x\of{1-c^2x}^3} \\
= \frac{c^4x^4}{\of{1-c^2x}^3}\\
\times\of{2 + 2c + cx - 4c^2x - 4c^3x + c^3x^2 + 2c^4x^2 + 
  2c^5x^2 - c^5x^3
}
\end{multline*}
which coincides with \eqref{eq:S312-2}.\qed

\bigskip\noindent
Now we turn to $\tau=\of{3,2,1}$.

\subsubsection{$\tau=\of{3,2,1}$, $r=2$, $\jv=(2)$}
By the same type of considerations as in our proof of
\eqref{eq:S321-1}, we observe that permutations with exactly one jump of depth 2 and exactly two occurrences of $\tau=\of{3,2,1}$ are in one--to--one
correspondence to lattice paths, where the
``path--segment'' around the jump looks {\em exactly\/} like
\epsffile{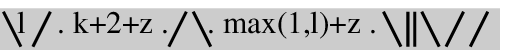},
at the end of which we are at height $(k+2)$ (the \tbase\ corresponding
to these permutations is $\of{4,3,1,2}$, see Table~\ref{tab:321} h).
Consider the generating function
of the ``lattice paths, starting from this path--segment up to the endpoint'':
For $l>1$, it is given as
\begin{equation*}
\sum_{k=0}^{\infty} \frac{x^{(k+4+l)/2}}{1-x} c^{k+3} x^{(k+2)/2} =
\frac{c^3x^{(l+6)/2}}{(1-c x)(1-x)}.
\end{equation*}
For $l=1$, we have
\begin{equation*}
\frac{c^3x^{9/2}}{(1-c x)(1-x)}.
\end{equation*}
For $l=0$, we have
\begin{equation*}
\frac{c^3x^{5}}{(1-c x)(1-x)}.
\end{equation*}
Thus, we may write down immediately the desired generating function
\begin{multline}
\label{eq:S321-2-2}
c\frac{c^3x^{5}}{(1-c x)(1-x)}+
c^3x^{3/2}\frac{c^3x^{9/2}}{(1-c x)(1-x)} \\+
\sum_{l=2}^{\infty}
	c^{l+2}x^{(l+1)/2}x^{1/2}\frac{c^3x^{(l+6)/2}}{(1-c x)(1-x)}
	\\ = \frac{c^4x^5\of{1-c x + c^2x + c^3x-c^3x^2}}{(1-x)(1-c x)^2}.
\end{multline}

\subsubsection{$\tau=\of{3,1,2}$, $r=2$, $\jv=(1)$} Recall the
considerations for the proof of \eqref{eq:S321-1}: 
There are exactly two ways to obtain precisely two occurrences from a
single jump of depth 1, namely
\begin{enumerate}
\item There is a second ``peak'' which is ``too near'' before the jump
	(see Lemma~\ref{lem:basic2}; the corresponding \tbase\ is $\of{3,4,2,1}$,
		see Table~\ref{tab:321} d),
\item There is a second down--step immediately after the jump,
	(\tbase\ $\of{4,2,3,1}$, see Table~\ref{tab:321} f).
\end{enumerate}
For the first case, note that the \tbase\ $\of{3,4,2,1}$, when
reflected at the upwards--sloping diagonal, is mapped to $\of{4,3,1,2}$:
Hence the generating function of this case is the same as in
\eqref{eq:S321-2-2}.
 
The generating function for the permutations corresponding to the
second case is $x \GFxtr{x}{\of{3,2,1}}{1}$: We simply have to ``insert''
one up--step before the jump and one down--step after the
jump in the path--segment we considered in the proof of \eqref{eq:S321-1}; i.e,
we must consider

\bigskip\centerline{\epsffile{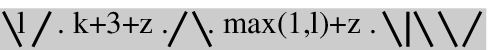}}\bigskip\noindent
here.

\subsubsection{$\tau=\of{3,2,1}$, $r=2$, $\jv=(1,1)$}
We want to determine the generating function for the set of lattice
paths which do contain precisely the two path--segments in the following way

\bigskip\centerline{\epsffile{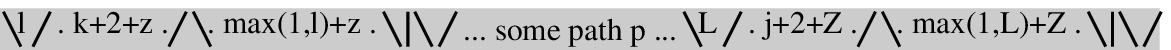},}\bigskip\noindent
where $p$ denotes some lattice path segment connecting the two
path--segments.

The corresponding \tbase s are $\of{3,2,1,6,5,4}$
and $(4,2,1,6,5,3)$; see Table~\ref{tab:321} a and e.

However, note that the ``connecting path--segment'' might be absent at all,
i.e., we face the situation 

\bigskip\centerline{\epsffile{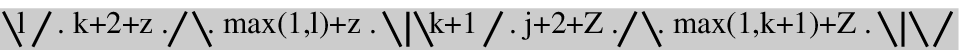}.}\bigskip\noindent

$\of{3,2,5,4,1}$, $(3,2,6,1,5,4)$,
$(4,2,6,1,5,3)$ and $(5,2,1,4,3)$ are the \tbase s corresponding
to this case, see Table~\ref{tab:321} b, c, g and i.

The same considerations as before finally lead to the generating function
\begin{multline}
\label{eq:S321-2-11}
\frac{c^3x^6\of{1-c x + c^2x}}{\of{1-x}^3\of{1-c x}^4}\\
\times\bigl(
	2-x-2c x+c^2x+cx^2+c^3x^2
	+c^4x^2-c^3x^3-2c^4x^3+c^4x^4
\bigr).
\end{multline}

\subsubsection{All together}
So, adding the respective generating functions \eqref{eq:S321-2-11},
divided by $x$, and 2 times \eqref{eq:S321-2-2}, divided by $x$, and
$\left(x \GFxtr{x}{\of{3,2,1}}{1}\right)$, we obtain
\begin{multline*}
\frac{c^3x^4}{\of{1-x}^3\of{1-c x}^4}\\
\times\bigl(
1 + 2\,c - 7\,c\,x - 5\,{c^2}\,x + 2\,{c^3}\,x + 
  2\,{c^4}\,x + 4\,c\,{x^2} + 16\,{c^2}\,{x^2} - 
  10\,{c^4}\,{x^2} - 4\,{c^5}\,{x^2} \\
  - c\,{x^3} - 
  10\,{c^2}\,{x^3} - 9\,{c^3}\,{x^3} + 
  15\,{c^4}\,{x^3} + 14\,{c^5}\,{x^3} + 
  2\,{c^6}\,{x^3} + 2\,{c^2}\,{x^4} + 
  6\,{c^3}\,{x^4} \\
  - 7\,{c^4}\,{x^4} - 
  16\,{c^5}\,{x^4} - 5\,{c^6}\,{x^4} - {c^3}\,{x^5} + 
  {c^4}\,{x^5} + 7\,{c^5}\,{x^5} + 4\,{c^6}\,{x^5} - 
  {c^5}\,{x^6} - {c^6}\,{x^6}
\bigr),
\end{multline*}
which after some simplification coincides with \eqref{eq:S321-2}.\qed


%% file: taubase312.tex
\begin{table}
\caption{$312$--bases and corresponding Dyck--paths with jumps for one
	and two occurrences}
\label{tab:312}
\vskip1em
\begin{tabular}{|l|c|c|}
\hline
Base $\rho$ & $\psi_{312}\of\rho$ & Reference in the text \\
\hline
\hline
\multicolumn{3}{|c|}{One occurrence:} \\
\hline
\hline
312 & \epsffile{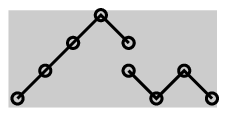} &   \\
\hline
\hline
\multicolumn{3}{|c|}{Two occurrences:} \\
\hline
\hline
312645 & \epsffile{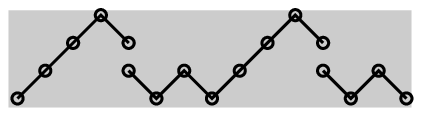} & a \\
\hline
31524 & \epsffile{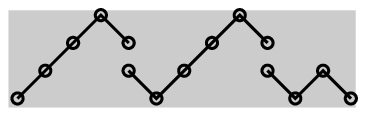} & b \\
\hline
316452 & \epsffile{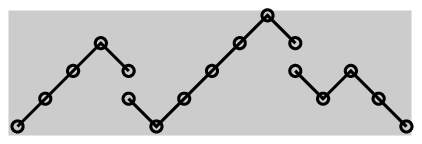} & c \\
\hline
3412 & \epsffile{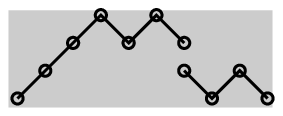} & d \\
\hline
4132 & \epsffile{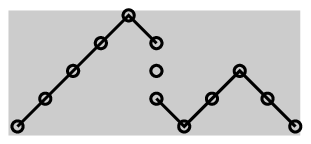} & e \\
\hline
4213 & \epsffile{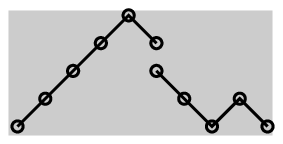} & f \\
\hline
423615 & \epsffile{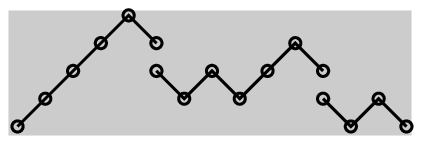} & g \\
\hline
4312 & \epsffile{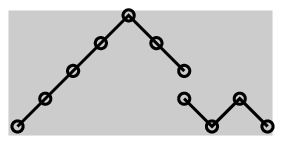} & h \\
\hline
\end{tabular}
\end{table}

%% file: taubase321.tex
\begin{table}
\caption{$321$--bases and corresponding Dyck--paths with jumps for one
	and two occurrences}
\label{tab:321}
\vskip1em
\begin{tabular}{|l|c|c|}
\hline
Base $\rho$ & $\psi_{321}\of\rho$ & Reference in the text \\
\hline
\hline
\multicolumn{3}{|c|}{One occurrence:} \\
\hline
\hline
321 & \epsffile{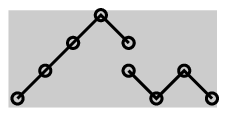} &  \\
\hline
\multicolumn{3}{|c|}{Two occurrences:} \\
\hline
\hline
321654 & \epsffile{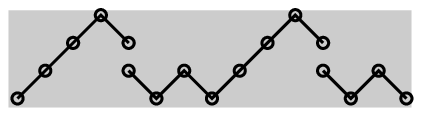} & a \\
\hline
32541 & \epsffile{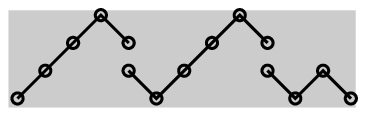} & b \\
\hline
326154 & \epsffile{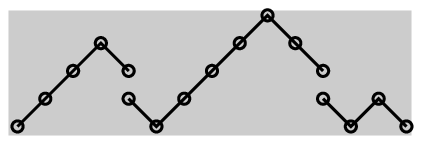} & c \\
\hline
3421 & \epsffile{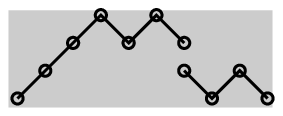} & d \\
\hline
421653 & \epsffile{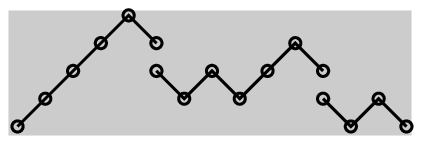} & e \\
\hline
4231 & \epsffile{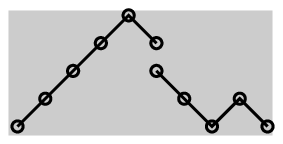} & f \\
\hline
426153 & \epsffile{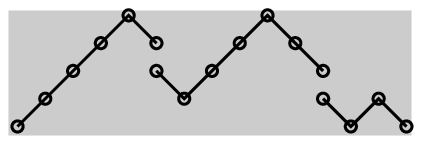} & g \\
\hline
4312 & \epsffile{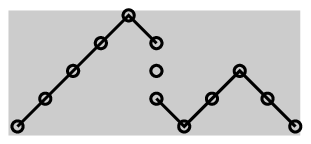} & h \\
\hline
52143 & \epsffile{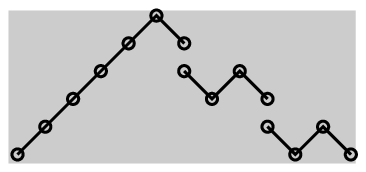} & i \\
\hline
\end{tabular}
\end{table}